\documentclass [12pt]{amsart}
\usepackage{times}
\usepackage{amssymb,latexsym,amsmath,amsfonts, eucal}
\usepackage[usenames,dvipsnames,svgnames,table]{xcolor}

\newcommand\hlight[1]{\tikz[overlay, remember picture,baseline=-\the\dimexpr\fontdimen22\textfont2\relax]\node[rectangle,fill=white!50,rounded corners,fill opacity = 0.2,draw,thick,text opacity =1] {$#1$};}

\usepackage[bookmarksnumbered, colorlinks, plainpages]{hyperref}

\newtheorem{theorem}{Theorem}[section]
\newtheorem{lemma}{Lemma}[section]
\newtheorem{remark}{Remark}[section]
\newtheorem{definition}{Definition}[section]
\newtheorem{example}{Example}[section]

\newtheorem{corollary}{Corollary}[section]
\numberwithin{equation}{section}
\usepackage{tikz}
\usetikzlibrary{shapes,arrows,positioning}
\usetikzlibrary{automata,arrows,positioning,calc}
\usepackage{enumitem}
\usepackage{subcaption}
\graphicspath{ {./images/} }

\title[ ]{Analysis of nonlinear fractal functions on PCF self-similar sets}
\author[Shah]{Aaryan Dharmesh Shah}
\address{Aaryan Dharmesh Shah, Department of Mathematics, National Institute of Technology Rourkela, Rourkela-769008, India}
\email{aaryanshah1911121@gmail.com}
\author[Jha]{Sangita Jha}
\address{Sangita Jha, Department of Mathematics, National Institute of Technology Rourkela, Rourkela-769008, India}
\email{jhasa@nitrkl.ac.in, \textbf{Corresponding Author-Sangita Jha}} 

\author[Mondal]{Anarul Islam Mondal}
\address{Anarul Islam Mondal, Department of Mathematics, National Institute of Technology Rourkela, Rourkela-769008, India}
\email{anarulmath96@gmail.com}

\subjclass[2010]{28A80, 35B41, 37L30, 47H10}

\keywords{PCF self-similar sets, fractal functions, Edelstein contraction, energy, normal derivatives, box-counting dimension}

\begin{document}
	
	
	\begin{abstract}
	
This article addresses three problems: (1) the construction of a general nonlinear fractal interpolation function on post-critically finite (PCF) self-similar sets; (2) the study of the energy and normal derivatives of uniform nonlinear fractal functions; and (3) the estimation of the box-counting dimensions of the proposed fractal functions on the Sierpi\'nski gasket and the von Koch curve.

We present a general framework for constructing attractors and fractal functions on PCF self-similar sets using Edelstein contractions. Further, we derive upper and lower bounds for the box-counting dimensions of the graphs of nonlinear interpolants. Finally, we provide several graphical and numerical examples to illustrate the construction and estimate the dimensions for different data sets.

	\end{abstract}
	
	\maketitle
	\tableofcontents

	%
	%
	%
	%
	\section{Introduction}
The present article studies several analytical properties of fractal interpolation functions (FIFs) on post-critically finite (PCF) self-similar sets. We first develop a general method for constructing fractal functions on PCF self-similar sets, and then study their analytical properties on some PCF sets such as Sierpi\'nski gasket ans von Koch curve. 

Fractal analysis is an important area of mathematical research that deals with complex structures showing self-similarity at different scales. One of its important applications is interpolation theory, which led to the development of fractal interpolation functions (FIFs). These functions were first introduced by Barnsley \cite{Barnsley1986} as a generalization of classical interpolation methods. Unlike classical interpolation functions, FIFs can describe irregular structures and may exhibit complex behaviour, including nowhere differentiability. Such functions are useful in modelling natural structures such as coastlines, mountain ranges, and biological patterns.

Since their introduction, several researchers have extended the theory of FIFs using different iterated function systems (IFSs). These include recurrent FIFs \cite{BarnsleyEltonHardin,BouDalla}, hidden-variable FIFs \cite{BEH}, piecewise FIF \cite{Chand1}, non-stationary FIFs \cite{Massopust2019}, \cite{MondalJha}, FIFs with variable
              scaling in subintervals \cite{Nasim}, and FIFs based on countable systems of data \cite{Secelean}.  A different approach to generalizing the classical IFS framework is to consider nonlinear IFS generated by suitable classes of weak contractions, including Meir--Keeler contractions \cite{Dumitru2009}, F-contractions \cite{Secelean2013}, Matkowski contractions \cite{Ri2017}, and Edelstein contractions \cite{Mihail, Pasupathi}.

Another approach is constructed from nonconstant harmonic functions goes back to the work of \c{C}elik, Ko\c{c}ak, and \"Ozdemir \cite{Celik}. They studied these curves as graphs of continuous functions on the Sierpi\'nski gasket (SG). In 2010, Ruan \cite{Ruan2010} extended the study of FIFs to a more general setting involving PCF self-similar fractals, a class of sets previously studied by Kigami \cite{Kigami}.

In 2013, de Amo, Carrillo, and S\'anchez  observed that nonconstant harmonic functions on ordinary intervals are always Lipschitz continuous, whereas this property does not necessarily hold for fractal sets such as the Sierpi\'nski gasket and the von Koch curve (KC). Using Barnsley's approach, they established the existence and uniqueness of invariant sets for IFSs on PCF self-similar structures under the assumption that the functions are Lipschitz continuous. However, this result does not cover IFSs on the SG constructed using nonconstant harmonic functions. This observations motivate us to develop a general framework for constructing fractal functions on PCF self-similar sets. 

The first main contribution of this article is the construction of nonlinear FIFs on PCF self-similar sets using Edelstein contractions. We prove the following (see Theorem \ref{graph}):

\begin{theorem}
		Let \(\mathcal{L} = (\mathcal{K}, \mathbb{N}_{N}, \{l_n\}_{n \in \mathbb{N}_{N}})\) be a PCF self-similar structure, and $\mathcal{I}$ be the  nonlinear IFS on $\mathcal{K}$. Then there exists a unique continuous function $f^*$ satisfying the interpolation condition (\ref{interpolation condition}) and whose graph is an invariant set of the IFS \(\mathcal{I}\). In particular, \[ f^*|_{V_{1}} = u \; \text{and}\;  f^*(l_{n}(x)) = s_{n}(f^*(x)) + h_{n}(x)\ \;\, \text{for all}\;\,\ n \in \mathbb{N}_{N},\; \ x \in \mathcal{K}.\]
		Moreover, the graph $ \Gamma f^*|_{\mathcal{K}}$ of \(f^*\) on \(\mathcal{K}\) satisfies
		\[\Gamma f^*|_{\mathcal{K}} = \bigcup\limits_{n \in \mathbb{N}_{N}} W_{n}(\Gamma f^*|_{l_{n}(\mathcal{K})}).\]
	\end{theorem}
Unlike the existing approaches that often rely on Lipschitz continuity and stronger contraction conditions, we use the H\"older continuity of nonconstant harmonic functions with respect to a suitable topologically equivalent metric. 
	
Further, we derive the bounds of the upper and lower box-counting dimensions of the graph of $f^*$ on SG and KC. To be specific, let $\mathcal K$ denote the SG and $G=\{(x,f^*(x)):x\in \mathcal K\}$ be the graph of $f^*$ defined in Theorem~\ref{graph}. The upper box-counting dimension on SG is the following (see Theorem \ref{SGUpperBox}):
\begin{theorem}
	Define
	$\mu=\sum_{n=1}^{3}L(s_n),$ where $L(s_n)$ denote the upper Lipschitz constants of the contraction maps $s_n$.
	Then
	\[
	\overline{\dim}_B(G)\leq
	\begin{cases}
		\frac{\log3}{\log2},
		&\mu\leq\frac{3}{2},\\
		1+\frac{\log\mu}{\log2},
		&\mu>\frac{3}{2}.
	\end{cases}
	\]
	In particular, if $\mu\leq3/2$, then $	\dim_B(G)=\frac{\log3}{\log2}.
	$
	\end{theorem}
	The following result provides a lower box-counting dimension on SG (see Theorem \ref{SGLowerBox}):
\begin{theorem}
 Define
	$\nu=\sum_{n=1}^{3}l(s_n)$, where $l(s_n)$ denote the lower Lipschitz constant of $s_n$.  Let
	\[
	S_1=\sum_{i=1}^{3}
	\operatorname{osc}_{l_i(\mathcal K)}f^*,
	\qquad
	H_1=\sum_{n=1}^{3}\sum_{i=1}^{3}
	\operatorname{osc}_{l_i(\mathcal K)}h_n.
	\]
	If $\nu>\frac32$ and
	$S_1>\frac{H_1}{\nu-\frac75}$, then
	$
	\underline{\dim}_B(G)
	\geq 1+\frac{\log\nu}{\log2}.
	$
\end{theorem}

The above theorems demonstrate that incorporating nonlinear contractions into the scaling functions leads to improved estimates for the box dimension of fractal interpolation functions (FIFs) on SG, compared with the existing results for linear fractal functions on SG (see \cite{Sahu,Manuj}). In particular, the authors of \cite{Sahu} and \cite{Manuj} established only upper bounds for the box dimension of FIFs (see \cite[Theorem 4.3]{Sahu} and \cite[Proposition 4]{Manuj}), using the H\"older exponent of the underlying functions on SG. Furthermore, to the best of our knowledge, no explicit dimension bounds have been established for fractal functions on KC.

Let $\mathcal{V}$ denote the KC and  $G=\{(x,f^*(x)):x\in\mathcal{V}\}$ be the graph of $f^*$ defined in Theorem~\ref{graph}. We obtain the upper and lower box-dimension of the graph of fractal function on KC (see Theorem \ref{KCUpperBox} and Theorem \ref{KCLowerBox})
\begin{theorem}	 
	 Define
	$\mu=\sum_{n=1}^{4}L(s_n)$. Then, \[
	\overline{\dim}_B(G)\leq
	\begin{cases}
		\frac{\log4}{\log3},
		&\mu\leq\frac43,\\[4pt]
		1+\frac{\log\mu}{\log3},
		&\mu>\frac43.
	\end{cases}
	\]
	In particular, if $\mu\leq4/3$, then
	$
	\dim_B(G)=\frac{\log4}{\log3}.
	$
\end{theorem}
\begin{theorem}
  Define
	$\nu=\sum_{n=1}^{4}l(s_n)$, and let
	\[
	S_1=\sum_{i=1}^{4}
	\operatorname{osc}_{l_i(\mathcal V)}f^*,
	\qquad
	H_1=\sum_{n=1}^{4}\sum_{i=1}^{4}
	\operatorname{osc}_{l_i(\mathcal V)}h_n.
	\]
	If $\nu>\frac43$ and
	$S_1>\frac{H_1}{\nu-1}$, then
	$
	\underline{\dim}_B(G)
	\geq1+\frac{\log\nu}{\log3}.
	$
\end{theorem}
The remainder of the paper is organized as follows. In Section \ref{sec-2}, we  provide a review of key definitions and results related to various fixed-point theorems, along with a discussion of PCF self-similar structures, energy forms, and harmonic functions.
	In Section \ref{sec-3}, we develop a framework for nonlinear fractal interpolation on PCF self-similar structures.  This includes illustrative examples of the constructed interpolants on SG and KC.
	In Section \ref{sec-4}, we estimate the box-counting dimension of the graphs of nonlinear FIFs defined on specific PCF fractals such as SG and KC.
Next we address our third problem in Section \ref{sec-5}. Here we explore properties such as energy of uniform nonlinear fractal functions on SG (see Theorem \ref{Thm-energySC}), KC (see Theorem \ref{Thm-energyKC}), and  normal derivatives on SG (see Theorem \ref{Thm-Normal}).

	\section{Preliminaries}\label{sec-2}
	In this section, we recall some well-known definitions and results that will be used later in this paper. We refer to \cite{Barnsley1993, Edelstein, MeirKeeler} for details.
	Throughout this paper, we use $\mathbb{N}^{0}=\mathbb{N} \cup \{0\},  \; \mathbb{N}_k$ as the set of the first $k$ natural numbers, and $\mathbb{N}_{k}^{0}=\mathbb{N}_{k} \cup \{0\}$.
	\begin{definition}
		Let $(X,d)$ be a metric space. A function $f:X \to X$ is called a Meir--Keeler contraction if for every $\epsilon > 0$ there exists $\delta > 0$ such that for 
		every $x, y \in X $,
			\[\epsilon \le d(x,y) < \epsilon + \delta \Rightarrow d(f(x),f(y)) < \epsilon.\]
	\end{definition}

	\begin{theorem}\label{MeirTheorem}
	Let $(X, d)$ be a complete metric space, and $f: X \to X$ be a mapping satisfying Meir--Keeler contraction. Then $f$ has a unique fixed point $\xi$. Moreover, for any $x \in X$,
	\[\lim\limits_{n \to \infty} f^{n}(x) = \xi.\]
\end{theorem}
	\begin{definition}
		Let $(X,d)$ be a metric space. Then $f:X \to X$ is called an Edelstein contraction if for every $x,y \in X$,
		\[x \neq y \implies d(f(x),f(y)) < d(x,y).\]
	\end{definition}
	\begin{remark}\label{MeirvsEdel}
		It is clear that every Meir--Keeler contraction is an Edelstein contraction. The converse holds in compact metric spaces (see \cite{MeirKeeler}).
	\end{remark}
\begin{definition}
	Let $(X,d)$ be a complete metric space and $f_i:X\to X$ be a contraction map for each $i=1,2,\dots,N $. Then $ \{X; f_i:i\in \mathbb{N}_N\}$ is known as a contractive iterated function system (IFS).
\end{definition}
	\begin{definition}(Lower and upper Lipschitz constants)
		Let $(X,d)$ be a metric space and let $f: X \rightarrow X$ be a function. Define $l(f)$ and $L(f)$ to be the lower and upper Lipschitz constants, respectively, where
		$$l(f)=\inf_{x,y\in X , x\neq y}\frac{d(f(x),f(y))}{d(x,y)},\; \; \text{and} \; \; L(f)=\sup_{x,y\in X,  x\neq y}\frac{d(f(x),f(y))}{d(x,y)}.$$
	\end{definition}
	\subsection{Post-critically finite self-similar structures}
	In this section, we recall some well-known facts of self-similar structure from \cite{Kigami}.
	\begin{definition}(Self-similar structure)
		Let \( \mathcal{K} \) be a compact metrizable topological space,  \( \Sigma = \mathbb{N}_{N}^{\mathbb{N}} \) is the one-sided shift space, and \( \{l_i : \mathcal{K} \to \mathcal{K}\}_{i \in \mathbb{N}_{N}} \) a collection of continuous injections on \( \mathcal{K} \). The tupple \( (\mathcal{K}, \mathbb{N}_{N}, \{l_i\}_{i \in \mathbb{N}_{N}}) \) is called a \emph{self-similar structure} if there exists a continuous surjection \( \pi: \Sigma \to \mathcal{K} \) such that for each \( i \in \mathbb{N}_{N},\; \; 	l_i \circ \pi = \pi \circ \sigma_i,$
		where the map \( \sigma_i: \Sigma \to \Sigma \) is defined by
		\[
		\sigma_i(w_1 w_2 w_3 \ldots) = i w_1 w_2 w_3 \ldots \quad \text{for } w_1 w_2 w_3 \ldots \in \Sigma.
		\]                                                            
		 For $\sigma=\sigma_1 \sigma_2 \sigma_3\ldots,\omega=\omega_1\omega_2\omega_3\ldots$, a metric on the one-sided shift space is defined by 
		 \[\delta_{\Sigma}(\sigma,\omega)=
		 \begin{cases} 0,\;\; \text{if}\;\; \sigma=\omega,\\
		 	  2^{-k}, \;\; \text{if}\;\; \sigma\neq \omega,
		 	 \end{cases}\]
	 	  where $k=\min\{n\in\mathbb{N}:\sigma_n\neq\omega_n\}$.
	\end{definition}
	
	\begin{definition}
		Let $\mathcal{L} = (\mathcal{K}, \mathbb{N}_{N}, \{l_i\}_{i \in \mathbb{N}_{N}})$ be a self-similar structure. 
		\begin{enumerate}
			\item The \textbf{critical set} $\mathcal{C}_{\mathcal{L}}$ is defined by
			\[
			\mathcal{C}_{\mathcal{L}} = \pi^{-1}(\mathcal{C}_{\mathcal{L}, \mathcal{K}}),\; \; \text{where}\;\; \mathcal{C}_{\mathcal{L}, \mathcal{K}} = \bigcup_{i, j \in \mathbb{N}_{N}, \, i \neq j} \big(l_i(\mathcal{K}) \cap l_j(\mathcal{K})\big).
			\]
			\item The \textbf{post-critical set} $\mathcal{P}_{\mathcal{L}}$ is defined by
			\[
			\mathcal{P}_{\mathcal{L}} = \bigcup_{n \geq 1} \sigma^n(\mathcal{C}_{\mathcal{L}}),
			\; \; \text{where}\; \; \sigma : \Sigma \to \Sigma\; \; \text{is the shift map given by}\; \; \]
			\[\sigma(w_1 w_2 w_3 \ldots) = w_2 w_3 \ldots \quad \text{for } w_1 w_2 w_3 \ldots \in \Sigma.
			\]
		\end{enumerate}
		Finally, define $V_0$ as the projection of the post-critical set
		$
		V_0 = \pi(\mathcal{P}_{\mathcal{L}}).$
	\end{definition}
	
	\begin{definition}
		Let \(\mathcal{L} = (\mathcal{K}, \mathbb{N}_{N}, \{l_i\}_{i \in \mathbb{N}_{N}})\) be a self-similar structure. Then \(\mathcal{L}\) is said to be post-critically finite (PCF for short) if and only if \(\mathcal{P}_{\mathcal{L}}\) is finite.
	\end{definition}
	Define the set of words of length \(n\) by \(W_{n} = \mathbb{N}_{N}^n =\{w_{1}w_{2}\dots w_{n}: n \in \mathbb{N}_{N}\}\), and \(W_{*} = \cup_{n \in \mathbb{N}^{0}} W_{n}\). Each word \(w \in W_{n}\) defines a continuous injection $l_{w}: \mathcal{K} \to \mathcal{K}$ by \(l_{w}:= l_{w_{1}}\circ l_{w_{2}}\circ\dots\circ l_{w_{n}}\). Set \( V_{n} = \cup_{i \in W_{n}} l_{i}(V_{0}), V_{*} = \displaystyle \cup_{n \in \mathbb{N}^{0}} V_{n}\). Then we have \(V_{n} \subset V_{n+1}\) and \( \mathcal{K} = {\overline{V}_{*}}\).
	\subsection{Energy, Laplacian, and harmonic function}
	In this section, we discuss the concepts Energy form, Laplace matrix, and harmonic function. 
	\begin{definition}(Laplace matrix)
		Let  $l(V_0):=\{f|f:V_0\rightarrow\mathbb{R}\}$. A symmetric linear operator $H:l(V_0)\xrightarrow{}l(V_0)$ is said to be a Laplace matrix on $V_0$, if
		\begin{enumerate}
			\item[(i)] $H_{pq}=H_{qp}\geq0$ for any $p\neq q$, where $H_{pq}:= H\chi_p(q),$
			\item[(ii)]  $(f,Hg):=\sum_{p\in V_0} f(p)\left(\sum_{q\in V_0} H_{pq}g(q)\right)\leq 0$ for any $f,g\in l(V_0)$, and
			\item[(iii)] $Hf=0$ if and only if $f$ is constant on $V_0$.
		\end{enumerate}
	\end{definition}
	\begin{definition}(Energy form)
		Let $H$ be a Laplace matrix, $S$ be a finite set and  $\textbf{r}:=\{r_i\}_{i\in S}$ be a set of positive numbers. We define a Dirichlet form or energy form $\mathcal{E}_m$ for $f,g\in V_m$ as follows
		$$\mathcal{E}_0(f,g)=-(f,Hg),$$
		$$\mathcal{E}_m(f,g)=\sum_{w\in S^m} r_w^{-1}\mathcal{E}_0(f|_{F_w(V_0)},g|_{F_w(V_0)}),$$ where $r_w=r_{w_1}\ldots r_{w_m}$ for $w=w_1w_2\ldots w_m\in S^m$.
		We set $\mathcal{E}_m(f):=\mathcal{E}_m(f,f).$
	\end{definition}    
	
	\begin{definition}(Regular harmonic structure)
		\label{reg}
		We say that $\mathcal{K}$ possesses a harmonic structure if there exists a pair $(H,\textbf{r})$ such that the following problem:
		$$\min\{\mathcal{E}_1(g): g|_{V_0}f\}=\mathcal{E}_0$$
		is solvable for all $f\in l(V_0)$ .	If $r_i\in \textbf{r}$ is less than $1$ for all $i\in S$, then the harmonic structure is said to be regular.
	\end{definition}
	If $\mathcal{K}$ possesses a regular harmonic structure, then it is clear that $\{\mathcal{E}_m\}_{m\geq0}$ is an increasing sequence. Define 
	$\mathcal{E}(f):=\lim_{m\rightarrow\infty}\mathcal{E}_m(f).$
	If the limit exists, we say that $f$ has finite energy.
	\begin{definition}(Effective resistance metric)
		For any $x,y\in V_m$, we define the effective resistance between $ x \; \text{and}\; y$ as
		  \[ R(x,y)^{-1}=
		\begin{cases}
			\min\{\mathcal{E}(f)| f(x)=0,f(y)=1\},\; \;  y\neq x,\\
			0,\; \; y=x.
		\end{cases}\] This defines a metric on $V_m$, hence on $V_*$, by continuity on $\mathcal{K}$, which is known as the resistance metric.
	\end{definition}
	From the definition of the resistance metric, it is seen that for any $f:\mathcal{K}\rightarrow\mathbb{R}$ having finite energy, the following holds:
	\[|f(x)-f(y)|^2\leq\mathcal{E}(f)R(x,y)\; \; \text{for all}\; \; x,y\in \mathcal{K}.\]
	
	\begin{definition}(Harmonic function)
		\label{harm}
		Let $f\in l(V_0)$. We define a harmonic function $h_f:\mathcal{K\rightarrow\mathbb{R}}$ to be the one that minimizes $\mathcal{E}_m(h_f)$ at each level $m\geq1$, where $h_f|_{V_0}=f$.
	\end{definition}
	Thus, from  Definitions \ref{reg} and \ref{harm}, it is easy to conclude that for a regular harmonic structure $\mathcal{K}$, any harmonic function $h$ is of finite energy and particularly:
	$$\mathcal{E}(h)=\mathcal{E}_0(h).$$
	
	\begin{definition}(Partition)
		Let $S$ be a finite set and let $\textbf{b}=\{b_i\}_{i\in S}$ be a family of numbers with $0 < b_i < 1$ for each $i\in S$. For $0 < \lambda < 1$, define		$$\Lambda_{\textbf{b}}(\lambda)=\{w=i_1i_2\ldots i_m\in \cup_{m\geq1} S^m: b_w\leq\lambda< b_{i_1i_2\ldots i_{m-1}}\}.$$
		We call $\Lambda_{\textbf{b}}(\lambda)$ the partition of $\Sigma$ with respect to $\textbf{b}$ and $\lambda$.
	\end{definition}
	The following is due to Hu and Wang \cite{Hu}, which establishes effective resistance estimates in terms of the underlying Euclidean metric on the fractal.
	\begin{theorem} \cite{Hu}
		\label{part}
		Let \(\mathcal{L} = (\mathcal{K}, \mathbb{N}_{N}, \{l_n\}_{n \in \mathbb{N}_{N}})\) be a PCF self-similar structure such that it possesses a regular harmonic structure $(H,\{r_i\}_{i=1}^n)$. Let $\textbf{b}=\{b_i\}_{i=1}^n$ exist with $ 0 < b_i < 1$ for every $i$ such that either
		\begin{enumerate}
			\item there exists a constant $ k_1> 0$ such that for any $ 0 < \lambda < 1:$
			\[\text{dist}(\mathcal{K}_w,\mathcal{K_\tau})\geq k_1\lambda \; \; \text{if}\; \; \mathcal{K}_w\cap\mathcal{K_\tau}=\phi \; \; \text{for}\;\; w,\tau\in\Lambda_b(\lambda),\] or
			\item there exist constants $ k_2,M> 0$ such that for any $0<\lambda<1$ and any $x_0\in\mathcal{K}$, any point $y\in B(x_0,k_2\lambda)$ can be connected to $x_0$ by a sequence of points$\{x_k\}_{k=0}^{n_0}$ in $\mathcal{K}$ with $1 \leq n_0 \leq M$, where $x_{n_0} = y$ and $x_{k-1},x_k\in \mathcal{K}_w,$ where $w\in \Lambda_\textbf{b}(\lambda)$.
		\end{enumerate}
		Then there exists $c>0$ such that
		\begin{equation*}
			R(x,y)\leq c|x-y|^{\alpha}\; \text{for all}\; x,y\in\mathcal{K},
		 \; \; \text{where}\;\; \alpha=\min_{1\leq i\leq n}\frac{\ln r_i}{\ln b_i}.
		 		 \end{equation*}
	\end{theorem}
	\section{A generalized fractal interpolation on PCF self-similar structures}\label{sec-3}
	Let \(\mathcal{K} \subseteq \mathbb{R}^{\tilde{n}}\), and let \(\mathcal{L} = (\mathcal{K}, \mathbb{N}_{N}, \{l_n\}_{n \in \mathbb{N}_{N}})\) denote a PCF self-similar structure. 
	Then, \[\mathcal{K} = \bigcup_{n \in \mathbb{N}_{N}} l_{n} (\mathcal{K}).\]
	We denote the elements of \(V_{0}\) by \(p_{1},p_{2},\dots,p_{L}\) and the elements of \(V_{1}\) by \(t_{1},t_{2},\dots,t_{M}\), under the assumption that \(p_{1} = t_{1},p_{2}=t_{2},\dots,p_{L}=t_{L}\).
	Let \(u: \mathbb{R}^{\tilde{n}} \supset V_{1} \to \mathbb{R}\) be any given function. Our aim is to construct a fractal interpolation function \(f\) on \(\mathcal{K}\) such that
	\begin{equation}\label{interpolation condition}
		f(t_{j}) = u(t_{j}) \ \ \text{for}\ j \in \mathbb{N}_{M}.
	\end{equation}
	For \(n \in \mathbb{N}_{N}\), let \(h_{n}: \mathcal{K} \to \mathbb{R}\) be a nonconstant harmonic function on \(\mathcal{K}\), \(s_{n}: \mathcal{J} \to \mathcal{J}\) be Edelstein contractions on the compact subset \(\mathcal{J}\) of \(\mathbb{R}\) and let \(F_{n}: \mathcal{K} \times \mathbb{R} \supset \mathcal{K} \times \mathcal{J} \to \mathcal{J}\) be \(N\) continuous maps defined by $F_{n}(x,y) = h_{n}(x) + s_{n}(y),$ 	and it satisfies $F_{n}(p_{k},u(t_{k})) = u(l_{n}(p_{k})).$
	That is, \begin{equation}\label{F cond}
		h_{n}(p_{k}) + s_{n}(u(t_{k})) = u(l_{n}(p_{k})), k=1,2,\ldots, L.
	\end{equation}
	
	We now define maps $W_{n}: \mathcal{K} \times \mathcal{J} \to \mathcal{K} \times \mathcal{J}, n \in \mathbb{N}_{N}$ as follows
	\[W_{n}(x,y) = (l_{n}(x), F_{n}(x,y))\ \text{for every}\ (x,y) \in \mathcal{K} \times \mathcal{J}.\]
	We call the IFS $\mathcal{I} = \{\mathcal{K} \times \mathcal{J};  W_{n}(x,y): n \in \mathbb{N}_{N}\}$ a nonlinear IFS on $\mathcal{K}$.
	Our main target is to find a unique continuous function on \(\mathcal{K}\) such that the function satisfies the above interpolation condition (\ref{interpolation condition}) and whose graph is the attractor of \(\mathcal{I}\). In the next theorem, we give the existence of an invariant set of the above nonlinear IFS on \(\mathcal{K}\). 
	
	\begin{theorem}\label{graph}
		Let \(\mathcal{L} = (\mathcal{K}, \mathbb{N}_{N}, \{l_n\}_{n \in \mathbb{N}_{N}})\) be a PCF self-similar structure, and $\mathcal{I}$ be the above defined nonlinear IFS on $\mathcal{K}$. Then there exists a unique continuous function $f^*$ satisfying the interpolation condition (\ref{interpolation condition}) and whose graph is an invariant set of the IFS \(\mathcal{I}\). In particular, \[ f^*|_{V_{1}} = u \; \text{and}\;  f^*(l_{n}(x)) = s_{n}(f^*(x)) + h_{n}(x)\ \; \text{for all}\; \ n \in \mathbb{N}_{N},\;\ x \in \mathcal{K}.\]
		Moreover, the graph of \(f^*\) on \(\mathcal{K}\), denoted by \(\Gamma f^*|_{\mathcal{K}}\), satisfies
		\[\Gamma f^*|_{\mathcal{K}} = \bigcup\limits_{n \in \mathbb{N}_{N}} W_{n}(\Gamma f^*|_{l_{n}(\mathcal{K})}).\]
	\end{theorem}
	
	
	\begin{proof}
		Let \(\mathcal{C}(\mathcal{K})\) be the set of all continuous functions on \(\mathcal{K}\). We consider the complete metric space
		\[\mathcal{C}^{*}(\mathcal{K}) = \{ g \in \mathcal{C}(\mathcal{K}) : g(p)= u(p)\, \text{ for all}\, p\in V_{0} \}\] endowed with the metric 
		\[d_{\infty} (f,g) = \sup\limits_{x \in \mathcal{K}} |f(x)-g(x)|.\]
		We define $ T : \mathcal{C}^{*}(\mathcal{K}) \to \mathcal{C}^{*}(\mathcal{K})$ by 
		\begin{equation} \label{RB operator eq}	
			(T g)(x) = h_{n}({l_{n}}^{-1}(x)) + s_{n}(g({l_{n}}^{-1}(x))),	\;(x \in l_{n}(\mathcal{K}),~~n \in \mathbb{N}_N).
		\end{equation}
		Let \(g \in \mathcal{C}^{*}(\mathcal{K})\). Since the functions \(h_n\), \(l_n^{-1}\), and \(g\) are continuous, it follows that \(Tg\) is continuous on each \(l_n(\mathcal{K})\); compatibility on overlaps follows from the join-up condition. Let \(p \in V_0\). Choose \(n \in \mathbb{N}_N\) such that \(p \in l_n(\mathcal{K})\), and set
		\(
		\tilde{p} := l_n^{-1}(p).
		\)
		Since
		\[
		V_0 \subset V_1 = \bigcup_{n\in \mathbb{N}_N} l_n(V_0),
		\]
		we have \(\tilde{p} \in V_0\). Therefore, using the join-up condition
		we obtain
		\begin{align*}
			Tg(p)
			&= h_n\big(l_n^{-1}(p)\big) + s_n\big(g(l_n^{-1}(p))\big)\\
			&= h_n(\tilde{p}) + s_n(u(\tilde{p}))= u(l_n(\tilde{p}))= u(p).
		\end{align*}
		Thus, \(Tg(p)=u(p)\) for every \(p \in V_0\). If \(p\) lies in more than one cell \(l_n(\mathcal{K})\), the same argument gives the same value \(u(p)\), so \(Tg\) is well-defined on \(V_0\). Thus, $Tg \in \mathcal{C}^{*}(\mathcal{K})$. 
		Since for each $n \in \mathbb{N}_{N},\ s_{n}: \mathcal{J} \to \mathcal{J}$ is an Edelstein contraction, it follows from Remark \ref{MeirvsEdel} that $s_{n}$ is also a Meir--Keeler contraction on $\mathcal{J}$. Let $\epsilon > 0$. Then, there exists $ \delta_{n} > 0$ such that for all $y_{1},y_{2} \in \mathcal{J}$,
		\begin{equation}\label{edelcond}
		 \epsilon \le |y_{1} - y_{2}| < \epsilon + \delta_{n}\ \ \ \implies |s_{n}(y_{1}) - s_{n}(y_{2})| < \epsilon.
		\end{equation}
		Let $g_{1}, g_{2} \in \mathcal{C}^{*}(\mathcal{K})$ be such that $\epsilon \le d_{\infty}(g_{1},g_{2}) < \epsilon + \delta,$ where $\delta = \min \{\delta_{n}: n \in \mathbb{N}_{N}\}$.
		Then \[|g_{1}(x)-g_{2}(x)| \le d_{\infty}(g_{1},g_{2}) < \epsilon + \delta \le \epsilon + \delta_{n}.\]
		Now, suppose 
		$\epsilon \le |g_{1}(x)-g_{2}(x)| \le d_{\infty}(g_{1},g_{2}) < \epsilon + \delta_{n}.$ Then  it follows from \eqref{edelcond}  that
		\[|s_{n}(g_{1}(x)) - s_{n}(g_{2}(x))| < \epsilon.\]
		Again, if $|g_{1}(x)-g_{2}(x)| < \epsilon \le d_{\infty}(g_{1},g_{2}) < \epsilon + \delta_{n},$ then
		\[|s_{n}(g_{1}(x)) - s_{n}(g_{2}(x))| < |g_{1}(x)-g_{2}(x)| < \epsilon.\]
		Thus in any case, $|s_{n}(g_{1}(x)) - s_{n}(g_{2}(x))| < \epsilon.$  Since the set $\mathcal{K}$ is compact and $g_{1}, g_{2}$ are continuous functions, \[\max\limits_{x \in \mathcal{K}} |s_{n}(g_{1}(x)) - s_{n}(g_{2}(x))| < \epsilon.\]
		Further,
		\begin{align*}
			d_{\infty}(Tg_{1},Tg_{2}) &= \max\limits_{x \in \mathcal{K}} |Tg_{1}(x) - Tg_{2}(x)|\\
			&= \max\limits_{n \in \mathbb{N}_{N}} \max\limits_{x \in l_{n}(\mathcal{K})} |s_{n}(g_{1}(l_{n}^{-1}(x))) - s_{n}(g_{2}(l_{n}^{-1}(x)))|\\
			&= \max\limits_{n \in \mathbb{N}_{N}} \max\limits_{x \in \mathcal{K}} |s_{n}(g_{1}(x)) - s_{n}(g_{2}(x))| < \epsilon.
		\end{align*}
		Therefore, $T$ is a Meir--Keeler contraction on $\mathcal{C}^{*}(\mathcal{K})$. By Theorem \ref{MeirTheorem}, $T$ therefore has a unique $f^{*}\in \mathcal{C}^{*}(\mathcal{K})$ such that $Tf^{*} = f^{*}$.
		
		Let $t' \in V_{1} \setminus V_{0}$ and let $t' \in l_{m}(\mathcal{K}) \cap l_{n}(\mathcal{K})$ for some $m,n \in \mathbb{N}_{N}$. Considering $t'$ as an element of $l_{n}(\mathcal{K})$, we obtain
		\begin{align*}
			f^{*}(t') &= h_{n}({l_{n}}^{-1}(t')) + s_{n}(f^{*}({l_{n}}^{-1}(t')))\\
			 & = h_{n}(p) + s_{n}(u(p))\; \text{for some}\ p \in V_{0}= u(t')\ \text{by using (\ref{F cond}}).
		\end{align*}
		Similarly, if we consider $t'$ as an element of $l_{m}(\mathcal{K})$, then we get $f^{*}(t') = u(t').$
		Hence, $f^{*}(t) = u(t)\ \text{for all}\ t \in V_{1}$.
		Furthermore, for each $n \in \mathbb{N}_{N}$
		\begin{align*}
			\Gamma f^*|_{l_{n}(\mathcal{K})} &= \{(x,h_{n}({l_{n}}^{-1}(x)) + s_{n}(f^{*}({l_{n}}^{-1}(x)))) : x \in l_{n}(\mathcal{K})\}\\
			&= \{(l_{n}(x),F_{n}(x,f^{*}(x)) : x \in \mathcal{K}\}\\
			&= \{W_{n}(x,f^{*}(x)): x \in \mathcal{K}\}= W_{n}(\{(x,f^{*}(x)): x \in \mathcal{K}\}).
		\end{align*}
		Thus, \begin{align*}
			\Gamma f^*|_{\mathcal{K}} &= \bigcup\limits_{n \in \mathbb{N}_{N}} W_{n}(\Gamma f^*|_{l_{n}(\mathcal{K})}).
		\end{align*}
		This completes the proof.
	\end{proof}
	
	\begin{remark}
		Using Barnsley's method, the authors of \cite{Amo} established a similar result by considering linear IFSs on PCF self-similar sets. Hence, our result can be regarded as a generalization of \cite[Theorem 3.6]{Amo}. 
	\end{remark}
	
	\begin{remark}
	Since every Banach contraction and every Matkowski contraction is, in particular, an Edelstein contraction; our result extends and unifies \cite[Theorem 2]{Celik}, where the mappings \( F_n \) are treated as Banach contractions in the second variable, as well as \cite[Theorem 3.3]{Ri2020}, in which the mappings are regarded as Matkowski contractions with respect to the second variable. 
	\end{remark}

 Our next result is dedicated to showing the uniqueness of the invariant set. But before that, we prove the following lemma.

\begin{figure}[h]
		
		\centering
		\begin{minipage}{0.48\textwidth}
			\centering
			\includegraphics[width=\linewidth]{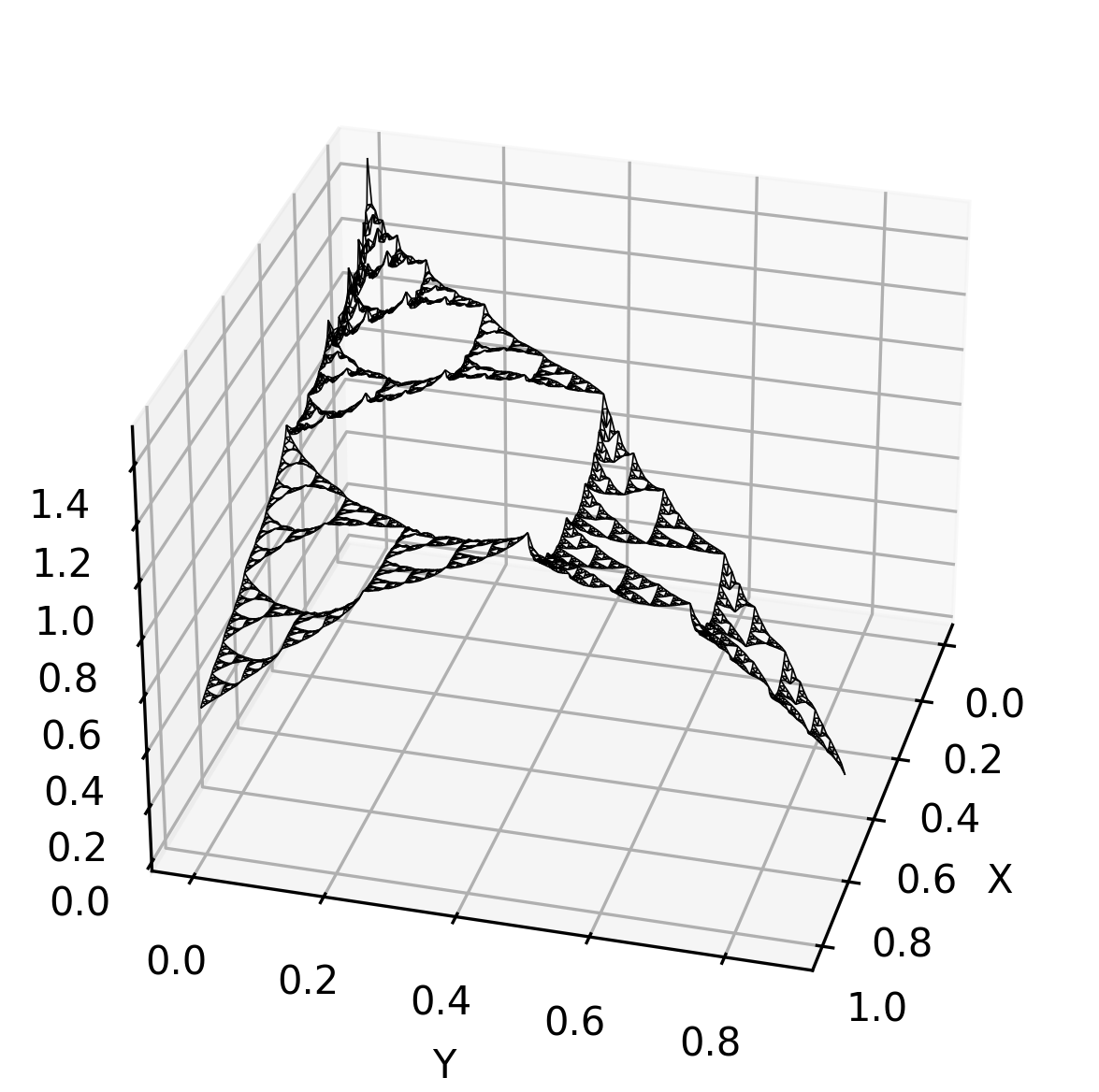}
			\captionof{figure}{A nonlinear FIF on the Sierpi\'nski gasket}
			\label{fig:nonlins}
		\end{minipage}
		\hfill
		\begin{minipage}{0.48\textwidth}
			
			\centering
			\includegraphics[width=\linewidth]{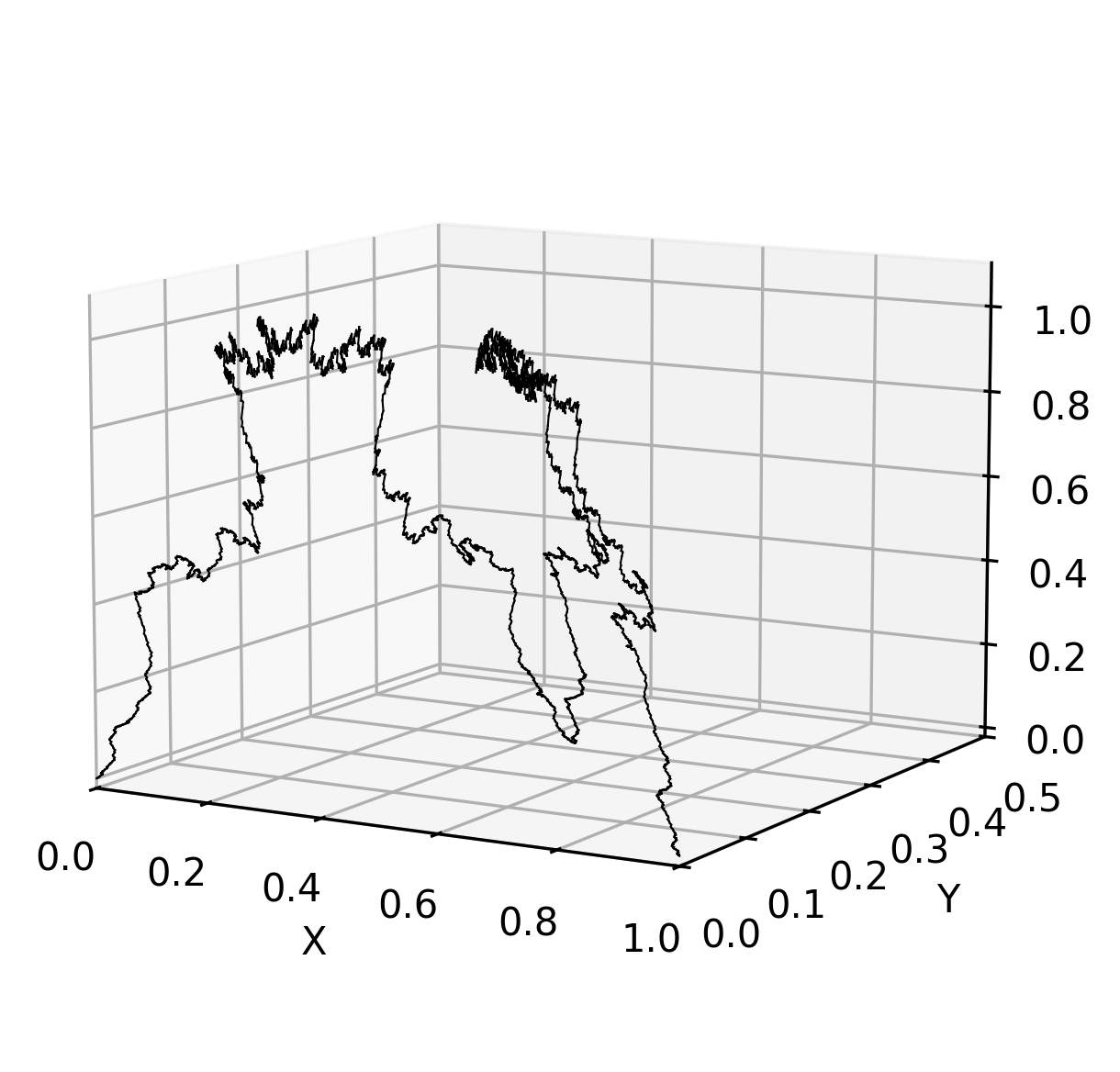}
			\captionof{figure}{A nonlinear FIF on the von Koch curve}
			\label{fig:nonlink}
		\end{minipage}
	\end{figure}

	\begin{lemma}\label{lemma3.1}
		Let \(\mathcal{L} = (\mathcal{K}, \mathbb{N}_{N}, \{l_n\}_{n \in \mathbb{N}_{N}})\) be a PCF self-similar structure equipped with a harmonic structure $(H,\{r_i\}_{i\in \mathbb{N}_N})$. Let there exist $\textbf{b}=\{b_i\}_{i\in \mathbb{N}_N}$ such that it satisfies one of the conditions given in Theorem \ref{part}. Let $h$ be a nonconstant harmonic function in $\mathcal{L}$. Then there exists $q>0 $ such that 
		$$|h(x)-h(x')|\leq q\cdot (|x_1-x_1'|^\alpha+|x_2-x_2'|^\alpha+\ldots+|x_{\tilde{n}}-x_{\tilde{n}}'|^\alpha)$$
		for all $x,x'\in \mathcal{K}\subset \mathbb{R}^{\tilde{n}}$, where  $\alpha=\frac{1}{2}\min_{1\leq i\leq n}\frac{\ln r_i}{\ln b_i}$.
	\end{lemma}
	\begin{proof}
		By the properties of the effective resistance metric and using the fact that a harmonic function has finite energy, we get
		$$|h(x)-h(x')|\leq\mathcal{E}(h)^{1/2}R(x,x')^{1/2}.$$
		Further, by Theorem \ref{part}, there is $c>0$ such that
		$$R(x,x')\leq c|x-x'|^{2\alpha}.$$
		Combining the above two inequalities, we get 
		\[|h(x)-h(x')|\leq (c\mathcal{E}(h))^{1/2}|x-x'|^{\alpha}\;\; \text{for all} \; \; x,x'\in \mathcal{K}.\]
		So $h$ is an $\alpha$-H\"older continuous function. Since $h$ can be any nonconstant harmonic function, it follows that $\alpha\leq1$.
		Therefore, by applying Jensen's inequality to the function $f(y)=y^\alpha$ and using the fact that $||x||_2\leq||x||_1$, we get
		$$|x-x'|^{\alpha}\leq|x_1-x_1'|^\alpha+|x_2-x_2'|^\alpha+\ldots+|x_{\tilde{n}}-x_{\tilde{n}}'|^\alpha .$$
	\end{proof}
	\begin{theorem}\label{Thm-3.2}
		Let \(\mathcal{L} = (\mathcal{K}, \mathbb{N}_{N}, \{l_n\}_{n \in \mathbb{N}_{N}})\) be a PCF self-similar structure equipped with a harmonic structure $(H,\{r_i\}_{i\in \mathbb{N}_N})$. Let there exist $\textbf{b}=\{b_i\}_{i\in \mathbb{N}_N}$ satisfying one of the conditions given in Theorem \ref{part}, and $\mathcal{I} = \{\mathcal{K} \times \mathcal{J}; W_{n}(x,y): n \in \mathbb{N}_{N}\}$ be the above defined nonlinear IFS on $\mathcal{K}$. Then there exists a metric $d_{\theta}$ on \(\mathcal{K} \times \mathcal{J} \subset \mathbb{R}^{\tilde{n}+1}\) topologically equivalent to the Euclidean metric $d$ with respect to which the maps $W_{n}$ are Edelstein contractions for each $n \in \mathbb{N}_{N}$.
	\end{theorem}
	\begin{proof}
		Define a metric $d_{\theta}: \mathcal{K} \times \mathcal{J}\to [0,\infty]$ by
		\begin{align*}
			d_{\theta}((x,y),(x',y')) &= d_{\theta}((x_{1},x_{2},\dots,x_{{\tilde{n}}},y),(x'_{1},x'_{2},\dots,x'_{{\tilde{n}}},y'))\\
			&= |x_{1}-x'_{1}|^\alpha+|x_{2}-x'_{2}|^\alpha+\dots+|x_{{\tilde{n}}}-x'_{{\tilde{n}}}|^\alpha+\theta|y-y'|,
		\end{align*}
		where $(x,y),(x',y') \in \mathcal{K} \times \mathcal{J},\alpha=\frac{1}{2}\min_{1\leq i\leq n}\frac{\ln r_i}{\ln b_i}, \theta \in \mathbb{R}^{+}$ specified below. Then, for all $(x,y),(x',y') \in \mathcal{K} \times \mathcal{J}$ with $(x,y) \neq (x',y')$, we have
		\allowdisplaybreaks
		\begin{align*}
			&d_{\theta}(W_{n}(x,y),W_{n}(x',y')) \\
			=& d_{\theta}((l_{n}(x),F_{n}(x,y)), (l_{n}(x'),F_{n}(x',y')))\\
			=& |l_{n,1}(x_{1})-l_{n,1}(x'_{1})|^\alpha+\dots+|l_{n,{\tilde{n}}}(x_{1})-l_{n,{\tilde{n}}}(x'_{1})|^\alpha+ \theta |F_{n}(x,y) - F_{n}(x',y')|\\
			\le& |a_{n}|^\alpha |x_{1}-x'_{1}|^\alpha +\dots+ |a_{n}|^\alpha |x_{{\tilde{n}}}-x'_{{\tilde{n}}}|^\alpha + \theta |h_{n}(x)-h_{n}(x')| + \theta |s_{n}(y) - s_{n}(y')|\\
			<& (|a_{n}|^\alpha + \theta q) ( |x_{1}-x'_{1}|^\alpha +\dots+ |x_{{\tilde{n}}}-x'_{{\tilde{n}}}|^\alpha)+ \theta |y-y'|, (\text{by using Lemma \ref{lemma3.1}})\\
			\le& \max\{|a_{n}|^\alpha + \theta q, 1\} d_{\theta}((x,y),(x',y'))\\
			\le& \max\{a^\alpha+\theta q, 1\} d_{\theta}((x,y),(x',y')),
		\end{align*}
		where $a = \max\limits_{n \in \mathbb{N}_{N}} \big\{|a_{n}|\big\}$ and $|a_{n}|$ are the contractivity factors of the maps $l_{n}$.\\
		Set $\theta = \dfrac{1-a^\alpha}{1+q}$. Then,
		\[d_{\theta}(W_{n}(x,y),W_{n}(x',y')) < d_{\theta}((x,y),(x',y'))\; \; \text{for all} \;\;(x,y),(x',y') \in \mathcal{K} \times \mathcal{J}\] with $ (x,y) \neq (x',y')$ and $ n \in \mathbb{N}_{N}.$ Hence, for each $n \in \mathbb{N}_{N}$, the map $W_{n}$ is Edelstein contraction with respect to the metric $d_{\theta}$.
	\end{proof}
	\begin{remark}
		Note that here we consider functions \(F_{n}\) that are H\"older continuous with respect to the first \({\tilde{n}}\) variables and Edelstein with respect to the last variable. A similar idea is used in \cite{Ri2017}, considering \(F_{n}\) to be Rakotch contractions with respect to the second variable. 
	\end{remark}
Next, we provide two illustrative examples.
	\begin{example}
		For each $n\in \{1,2,3\}$, let $l_n: \mathbb{R}^2\to \mathbb{R}^2$ be given by $l_n=\frac{1}{2}(x-q_n)+q_n$, where $q_1=(0,0),q_2=(1,0), q_3=(0.5,\frac{\sqrt{3}}{2})$ and  $s_n:[0,2]\to [0,2]$ be defined by
		$ s_1(y)=\frac{1}{4}y^2,\; s_2(y)=\frac{1}{2}y,\; s_3(y)=\frac{1}{2+y}.$
				Consider the interpolation data \[
		\begin{aligned}
			&\left\{ (0, 0, 0.5), \left(0.5, 0, 1\right), (1, 0, 1.5),\left(0.25, \frac{\sqrt{3}}{4}, 1\right), \left(0.75, \frac{\sqrt{3}}{4}, 1\right), \left(\frac{0.5,\sqrt{3}}{2}, 0\right) \right\}.
		\end{aligned}
		\] The data obeys the following functional equation
		$$f^*(l_{n}(x)) = s_{n}(f^*(x)) + h_{n}(x)\ \; \text{for  all}\; \ n \in \{1,2,3\} \;\text{and}\ x \in \mathcal{K},$$
		where $h_n$ is the harmonic function determined by the given interpolation data and the choice of $s_n$, and $ \mathcal{K}$ is the SG. 
		The figure of the nonlinear FIF on the SG is given in Fig \ref{fig:nonlins}.
	\end{example}
	\begin{example}
		For $n\in \{1,2,3,4\}$, let $l_n:\mathbb{R}^2\to \mathbb{R}^2$ be defined by
		$${l_1(x,y)=\frac{1}{3}(x,y)},$$
		$${l_2(x,y)=\frac{1}{3}\left(\frac{x-\sqrt{3}y}{2},\frac{\sqrt{3}x+y}{2}\right)+\left(\frac{1}{3},0\right)},$$
		$$l_3(x,y)=\frac{1}{3}\left(\frac{x+\sqrt{3}y}{2},\frac{y-\sqrt{3}x}{2}\right)+\left(\frac{1}{2},\frac{1}{2\sqrt{3}}\right),$$
		$$l_4(x,y)=\frac{1}{3}(x,y)+\left(\frac{2}{3},{0}\right),$$
		and $q_1=(0,0),q_2=(1,0)$.
		Consider the interpolation data \[
		\begin{aligned}
			&\left\{ (0, 0, 0), \left(\frac{1}{3}, 0, 1\right), \left(\frac{1}{2},\frac{1}{2\sqrt{3}} , 0\right),\left(\frac{2}{3}, 0, 1\right), \left(1, 0, 0\right) \right\}
		\end{aligned}
		\] where the interpolant satisfies the following functional equation:		
		$$f^*(l_{n}(x)) = s_{n}(f^*(x)) + h_{n}(x)\ \,  \text{for all}\, \ n \in \{1,2,3\}\ \text{and}\ x \in \mathcal{K},$$
		where $\mathcal{K}$ is the KC, and  $s_n:[0,\frac{6}{5}]\to [0,\frac{6}{5}]$ are given by
		\begin{align*}
			s_1(y)=\frac{1}{4}y^2,\qquad s_2(y)=\frac{1}{4}y^2,\qquad s_3(y)=\frac{1}{4}y^2,\qquad s_4(y)=\frac{1}{4}y^2.
		\end{align*}
		Here, $h_n$ is the harmonic function determined by the given interpolation data and the choice of $s_n$.
		The figure of the nonlinear FIF on the KC is given in Fig \ref{fig:nonlink}.
	\end{example}
	
	\section{Dimensional analysis on certain fractals}\label{sec-4}
	In this section, our main goal is to estimate the box dimension of the graph of nonlinear FIFs generated on PCF self-similar structures. For this, we have taken two PCF self-similar sets, namely SG and KC. Throughout this section, $\mathcal{K}$ denotes the SG, and $\mathcal{V}$ denotes KC. To our knowledge, there are currently no studies addressing both the upper and lower dimensions of nonlinear FIFs on the SG or KC.
	\begin{definition}(Box-counting dimension) {
			   Let $F\subset\mathbb{R}^{\tilde{n}}$ be a nonempty set, and let $N_\delta(F)$ denote the minimum number of cubes of side length $\delta$ that cover $F$.	The lower box-counting dimension is defined by
			$$ \underline{\dim}_B(F)=\liminf_{\delta \rightarrow0} \frac{\log(N_{\delta}(F))}{-\log\delta},$$
			and the upper box-counting dimension is defined by
			$$ \overline\dim_B(F)=\limsup_{\delta \rightarrow0} \frac{\log(N_{\delta}(F))}{-\log\delta}.$$
		}
	\end{definition}
	If both of these values are equal, we define the common value as the box-counting dimension of $F$.
	\begin{definition}(Oscillation)
	Let $\mathcal{J}\subset \mathbb{R}^n$ be a nonempty set and let $f:\mathcal{J}\xrightarrow{}\mathbb{R}$ be a function. The oscillation of the function $f$ on its domain $\mathcal{J}$ is defined by
	\[
	\operatorname{osc}_{\mathcal J}f
	:=\sup\bigl\{|f(x)-f(y)|:x,y\in\mathcal J\bigr\}.
	\]
	\end{definition}
	\subsection{Box-counting dimensions of nonlinear FIFs on SG}
	\begin{lemma}
		\label{SGcover}
		Let $f^*:\mathcal{K}\to\mathbb{R}$
		be continuous. Write
		\[
		G=\{(x,f^*(x)):x\in\mathcal{K}\},
		\qquad
		G_w=\{(x,f^*(x)):x\in l_w(\mathcal{K})\}.
		\]
		Let $N(k)$ and $N_w(k)$ denote the minimum numbers of
		cubes of side length $2^{-k}$ required to cover $G$
		and $G_w$, respectively. Then
		\[
		\frac1{21}\sum_{|w|=k}N_w(k)
		\leq N(k)
		\leq\sum_{|w|=k}N_w(k).
		\]
	\end{lemma}
	\begin{proof}
		Since $G=\bigcup_{|w|=k}G_w$, combining minimal covers
		of the cell graphs gives
		\[
		N(k)\leq\sum_{|w|=k}N_w(k).
		\]
		For the reverse inequality, put $\varepsilon=2^{-k}$
		and choose a minimal cover $\{Q_j\}_{j=1}^{N(k)}$ of
		$G$ by cubes of side length $\varepsilon$. Define
		\[
		I_w=\{j:Q_j\cap G_w\neq\varnothing\}.
		\]
		These cubes cover $G_w$, so $N_w(k)\leq\#I_w$.
		Therefore, we have
		\[
		\sum_{|w|=k}N_w(k)
		\leq\sum_{j=1}^{N(k)}
		\#\{w:|w|=k,\ Q_j\cap G_w\neq\varnothing\}.
		\]
		
		Let $T$ be the closed equilateral triangle bounded
		by the three vertices of $\mathcal{K}$.
		The triangles $T_w=l_w(T)$, $|w|=k$, have side
		length $\varepsilon$ and pairwise disjoint interiors.
		If $Q_j$ meets $G_w$, its horizontal projection meets
		$T_w$. Since $\operatorname{diam}(T_w)=\varepsilon$,
		every such triangle lies in the square obtained by
		extending the horizontal projection of $Q_j$ by
		$\varepsilon$ in each coordinate direction. This enlarged square has an area of $9\varepsilon^2$,
		whereas each $T_w$ has an area
		 of $\sqrt{3}\varepsilon^2/4$. Hence, the number of these
		triangles are at most
		\[
		\frac{9\varepsilon^2}
		{\sqrt{3}\varepsilon^2/4}
		=12\sqrt{3}<21.
		\]
		Therefore,
		\[
		\sum_{|w|=k}N_w(k)\leq21N(k).
		\]
	\end{proof}
	 \begin{lemma}\label{SGharmosc}
		Let $h$ be a harmonic function on the Sierpi\'nski gasket $\mathcal{K}$, then we have
		\[
		\sum_{i=1}^{3}\operatorname{osc}_{l_i(\mathcal{K})} h \leq\frac{7}{5}\operatorname{osc}_{\mathcal{K}}h.
		\]
	\end{lemma}
	\begin{proof}
		If $h$ is constant, then we are done. Otherwise, we can normalize the boundary values of function such that $\operatorname{osc}_{\mathcal{K}}h=1$ by taking $
		h(q_1)=0,\qquad h(q_2)=t,\qquad h(q_3)=1,
		 0\le t\le1.
		$
		The harmonic extension rule gives
		\[
		h(l_1(q_2))=h(l_2(q_1))=\frac{2t+1}{5},\;
		h(l_1(q_3))=h(l_3(q_1))=\frac{t+2}{5},\;
		h(l_2(q_3))=h(l_3(q_2))=\frac{2t+2}{5}.
		\]
		Now, the extension rule implies that any harmonic function must attain its maximum at the boundary values. Hence,
		\[
		\operatorname{osc}_{l_1(\mathcal K)}h=\frac{t+2}{5},\;	\operatorname{osc}_{l_3(\mathcal K)}h=\frac{3-t}{5}.
		\]
		Moreover,
		\begin{align*}
			\operatorname{osc}_{l_2(\mathcal K)}h
			&=\max\left\{t,\frac{2t+1}{5},\frac{2t+2}{5}\right\}
			-\min\left\{t,\frac{2t+1}{5},\frac{2t+2}{5}\right\}\\
			&=\frac{1}{5}\max\{1,\,2-3t,\,3t-1\}
			\le\frac{2}{5}.
		\end{align*}
		Therefore,
		\[
		\sum_{i=1}^{3}\operatorname{osc}_{l_i(\mathcal K)}h
		\le\frac{t+2}{5}+\frac25+\frac{3-t}{5}
		=\frac{7}{5}=\frac{7}{5}\operatorname{osc}_{\mathcal{K}}h.
		\]
	\end{proof}
	Now, we give upper and lower bounds on the box dimension of nonlinear functions constructed on $\mathcal{K}$.
	\begin{theorem}\label{SGUpperBox}
	Let $G=\{(x,f^*(x)):x\in\mathcal{K}\}$ be the graph
	of $f^*$ defined in Theorem~\ref{graph}. Define
	$\mu=\sum_{n=1}^{3}L(s_n).$
	Then
	\[
	\overline{\dim}_B(G)\leq
	\begin{cases}
		\frac{\log3}{\log2},
		&\mu\leq\frac{3}{2},\\
		1+\frac{\log\mu}{\log2},
		&\mu>\frac{3}{2}.
	\end{cases}
	\]
	In particular, if $\mu\leq3/2$, then
	\[
	\dim_B(G)=\frac{\log3}{\log2}.
	\]
	\end{theorem}
	\begin{proof}
		Define for $m\geq1$
		\[
		S_m=\sum_{|w|=m}
		\operatorname{osc}_{l_w(\mathcal K)}f^*,
		\qquad
		H_m=\sum_{n=1}^{3}\sum_{|w|=m}
		\operatorname{osc}_{l_w(\mathcal K)}h_n.
		\]
		For any word $w$ of length $m$, $h_n\circ l_w$
		is also harmonic. Applying Lemma~\ref{SGharmosc}, we have
		\begin{align*}
			H_{m+1}
			&=\sum_{n=1}^{3}\sum_{|w|=m}\sum_{i=1}^{3}
			\operatorname{osc}_{l_{wi}(\mathcal K)}h_n\\
			&\leq\frac75
			\sum_{n=1}^{3}\sum_{|w|=m}
			\operatorname{osc}_{l_w(\mathcal K)}h_n
			=\frac75H_m.
		\end{align*}
		Applying this relation recursively, we obtain
		\[
		H_m\leq H_1\left(\frac75\right)^{m-1}.
		\]
		For every word $w$ of length $m$, the functional
		equation yields
		\begin{align*}
			\operatorname{osc}_{l_{nw}(\mathcal K)}f^*
			&=\operatorname{osc}_{l_w(\mathcal K)}
			\bigl(s_n\circ f^*+h_n\bigr)\\
			&\leq
			L(s_n)\operatorname{osc}_{l_w(\mathcal K)}f^*
			+\operatorname{osc}_{l_w(\mathcal K)}h_n.
		\end{align*}
		Summing over $n$ and $w$, we get
		\begin{align*}
			S_{m+1}
			&=\sum_{n=1}^{3}\sum_{w\in W_m}
			\operatorname{osc}_{l_{nw}(\mathcal K)}f^*\\
			&\leq\sum_{n=1}^{3}\sum_{w\in W_m}
			\left[
			L(s_n)\operatorname{osc}_{l_w(\mathcal K)}f^*
			+\operatorname{osc}_{l_w(\mathcal K)}h_n
			\right]\\
			&=\left(\sum_{n=1}^{3}L(s_n)\right)S_m+H_m\\
			&=\mu S_m+H_m\\
			&\leq\mu S_m+
			H_1\left(\frac75\right)^{m-1}.
		\end{align*}
		Let $N(m)$ denote the minimum number of cubes of
		side length $2^{-m}$ required to cover $G$, and let
		$N_w(m)$ denote the corresponding covering number
		of the graph over $l_w(\mathcal K)$. Each level-$m$ cell has a diameter of $2^{-m}$. Therefore, the graph over this cell is contained
		in a rectangular prism with horizontal side
		lengths at most $2^{-m}$ and height
		$\operatorname{osc}_{l_w(\mathcal K)}f^*$.
		Hence, there exists $C_0>0$, independent
		of $m$ and $w$, such that
		\[
		N_w(m)\leq C_0\left(
		1+2^m\operatorname{osc}_{l_w(\mathcal K)}f^*
		\right)
		\]
		From Lemma \ref{SGcover}, we obtain
		\begin{align*}
			N(m)
			&\leq\sum_{w\in W_m}N_w(m)\\
			&\leq C_0\sum_{w\in W_m}
			\left(
			1+2^m\operatorname{osc}_{l_w(\mathcal K)}f^*
			\right)\\
			&=C_0\left(3^m+2^mS_m\right).
		\end{align*}
		
		\noindent
		\textbf{Case (i): $\mu\leq3/2$.}
		Since $S_m\geq0$, we have
		\[
		S_{m+1}\leq\frac32S_m+
		H_1\left(\frac75\right)^{m-1}.
		\]
		Iterating this inequality, we obtain
		\begin{align*}
			S_m
			&\leq
			\left(\frac32\right)^{m-1}S_1
			+H_1\sum_{j=1}^{m-1}
			\left(\frac32\right)^{m-1-j}
			\left(\frac75\right)^{j-1}\\
			&=
			\left(\frac32\right)^{m-1}S_1
			+
			H_1
			\frac{
				\left(\frac32\right)^{m-1}
				-\left(\frac75\right)^{m-1}
			}{
				\frac32-\frac75
			}\\
			&\leq
			\left(\frac32\right)^{m-1}
			\left(
			S_1+\frac{H_1}{\frac32-\frac75}
			\right).
		\end{align*}
		Hence, there exists $C_1>0$, independent of $m$,
		such that
		\[
		S_m\leq C_1\left(\frac32\right)^{m-1}.
		\]
		Substituting this into the covering estimate gives
		\begin{align*}
			N(m)
			&\leq C_0\left[
			3^m+C_1\,2^m
			\left(\frac32\right)^{m-1}
			\right]\\
			&=C_0\left(1+\frac{2C_1}{3}\right)3^m.
		\end{align*}
		Therefore,
		\[
		\overline{\dim}_B(G)
		=\limsup_{m\to\infty}
		\frac{\log N(m)}{m\log2}
		\leq\frac{\log3}{\log2}.
		\]
		
		\noindent
		\textbf{Case (ii): $\mu>3/2$.}
		Iterating the oscillation inequality gives
		\begin{align*}
			S_m
			&\leq\mu^{m-1}S_1
			+H_1\sum_{j=1}^{m-1}
			\mu^{m-1-j}\left(\frac75\right)^{j-1}\\
			&=\mu^{m-1}S_1
			+H_1
			\frac{\mu^{m-1}-\left(\frac75\right)^{m-1}}
			{\mu-\frac75}\\
			&\leq\mu^{m-1}
			\left(S_1+\frac{H_1}{\mu-\frac75}\right).
		\end{align*}
		Let
		\[
		C_2=S_1+\frac{H_1}{\mu-\frac75}.
		\]
		Since $2\mu>3$, we obtain
		\begin{align*}
			N(m)
			&\leq C_0\left(3^m+2^mS_m\right)\\
			&\leq C_0\left(
			3^m+\frac{C_2}{\mu}(2\mu)^m
			\right)\\
			&\leq C_0\left(1+\frac{C_2}{\mu}\right)(2\mu)^m.
		\end{align*}
		Therefore, we get
		\[
		\overline{\dim}_B(G)
		=\limsup_{m\to\infty}
		\frac{\log N(m)}{m\log2}
		\leq1+\frac{\log\mu}{\log2}.
		\]
		
		On the other hand, projection onto $\mathcal K$ is Lipschitz
		and surjective, so
		\[
		\underline{\dim}_B(G)
		\geq\dim_B(\mathcal K)
		=\frac{\log3}{\log2}.
		\]
		Thus, if $\mu\leq3/2$, the box dimension exists
		and satisfies
		\[
		\dim_B(G)=\frac{\log3}{\log2}.
		\]
	\end{proof}

\begin{theorem}\label{SGLowerBox}
	Let $G=\{(x,f^*(x)):x\in\mathcal{K}\}$ be the graph
	of $f^*$ defined in Theorem~\ref{graph}. Define
	$\nu=\sum_{n=1}^{3}l(s_n)$, and let
	\[
	S_1=\sum_{i=1}^{3}
	\operatorname{osc}_{l_i(\mathcal K)}f^*,
	\qquad
	H_1=\sum_{n=1}^{3}\sum_{i=1}^{3}
	\operatorname{osc}_{l_i(\mathcal K)}h_n.
	\]
	If $\nu>\frac32$ and
	$S_1>\frac{H_1}{\nu-\frac75}$, then
	\[
	\underline{\dim}_B(G)
	\geq 1+\frac{\log\nu}{\log2}.
	\]
\end{theorem}
	\begin{proof}
		Define for $m\geq1$
		\[
		S_m=\sum_{|w|=m}\operatorname{osc}_{l_w(\mathcal K)}f^*,
		\qquad
		H_m=\sum_{n=1}^{3}\sum_{|w|=m}
		\operatorname{osc}_{l_w(\mathcal K)}h_n.
		\]
		For any word $w$ of length $m$, $h_n \circ l_w$ is also harmonic. Applying Lemma \ref{SGharmosc}, we have
		
		\begin{align*}
			H_{m+1}
			&=\sum_{n=1}^{3}\sum_{|w|=m}\sum_{i=1}^{3}
			\operatorname{osc}_{l_{wi}(\mathcal K)}h_n\\
			&\le\frac{7}{5}
			\sum_{n=1}^{3}\sum_{|w|=m}
			\operatorname{osc}_{l_w(\mathcal K)}h_n
			=\frac{7}{5}H_m.
		\end{align*}
		Applying this relation recursively, we have
		\[
		H_m\leq H_1\left(\frac{7}{5}\right)^{m-1}.
		\]
		For every word $w$ of length $m$, the functional equation yields
		\begin{align*}
			\operatorname{osc}_{l_{nw}(\mathcal K)}f^*
			&=\operatorname{osc}_{l_w(\mathcal K)}
			\bigl(s_n\circ f^*+h_n\bigr)\\
			&\ge
			\ell(s_n)\operatorname{osc}_{l_w(\mathcal K)}f^*
			-\operatorname{osc}_{l_w(\mathcal K)}h_n.
		\end{align*}
		Summing over $n$ and $w$, we get
	
		\begin{align*}
			S_{m+1}
			&=\sum_{n=1}^{3}\sum_{w\in W_m}
			\operatorname{osc}_{l_{nw}(\mathcal K)}f^*\\
			&\ge
			\sum_{n=1}^{3}\sum_{w\in W_m}
			\left[
			l(s_n)\operatorname{osc}_{l_w(\mathcal K)}f^*
			-
			\operatorname{osc}_{l_w(\mathcal K)}h_n
			\right]\\
			&=
			\sum_{n=1}^{3}l(s_n)
			\left[
			\sum_{w\in W_m}
			\operatorname{osc}_{l_w(\mathcal K)}f^*
			\right]
			-
			\sum_{n=1}^{3}\sum_{w\in W_m}
			\operatorname{osc}_{l_w(\mathcal K)}h_n\\
			&=
			\left(\sum_{n=1}^{3}l(s_n)\right)S_m-H_m\\
			&=\nu S_m-H_m\\
			&\ge\nu S_m-H_1\left(\frac{7}{5}\right)^{m-1}.
		\end{align*}
		Iterating the above inequality over $m$, we obtain
		\begin{align*}
			S_m
			&\ge\nu^{m-1}S_1
			-H_1\sum_{j=1}^{m-1}\nu^{m-1-j}\left(\frac{7}{5}\right)^{j-1}\\
			&=\nu^{m-1}
			\left(S_1-\frac{H_1}{\nu-7/5}\right)
			+\frac{H_1\left(\frac{7}{5}\right)^{m-1}}{\nu-7/5}\\
			&=C_1\nu^{m-1} + C_2 \left(\frac{7}{5}\right)^{m-1},
		\end{align*}
		where $C_1 = \left(S_1-\frac{H_1}{\nu-7/5}\right)$ and $C_2=\left(\frac{H_1}{\nu-7/5}\right)$.
		Since we have assumed that $C_1$ is positive, we have
		\[ 
		S_m\ge C_1\nu^{m-1}.
		\]
		For each $w\in W_m$, let $N_w(m)$ denote the minimum
		number of cubes of side length $2^{-m}$ required to
		cover the graph of $f^*$ over $l_w(\mathcal K)$.
		Since $l_w(\mathcal K)$ is connected and $f^*$ is
		continuous, the vertical projection of this graph
		is an interval of length
		$\operatorname{osc}_{l_w(\mathcal K)}f^*$.
		Each covering cube has a vertical projection of
		length $2^{-m}$. Therefore,
		\[
		N_w(m)\geq
		2^m\operatorname{osc}_{l_w(\mathcal K)}f^*.
		\]
		Summing over $w\in W_m$ and applying
		Lemma \ref{SGcover}, we obtain
		\begin{align*}
			N(m)
			&\geq C\sum_{w\in W_m}N_w(m)\\
			&\geq C\,2^m\sum_{w\in W_m}
			\operatorname{osc}_{l_w(\mathcal K)}f^*\\
			&=C\,2^mS_m,
		\end{align*}
		where $C>0$ is independent of $m$. Putting in the equation of $S_m$, we have
		\[
		N(m)\ge C (2\nu)^m,\; \; \text{where}\;\; C= \frac{\cdot C_1}{\nu}.\]
		Therefore, we obtain
		\[
		\dim_B(G)= \liminf_{m\xrightarrow{}\infty}\frac{\log N(m)}{\log 2^m}\ge \liminf_{m\xrightarrow{}\infty}\frac{\log C(2\nu)^m}{\log 2^m} =1 + \frac{\log \nu}{\log 2}.
		\]
		
	\end{proof}
	
	\subsection{Examples}
\begin{figure}[h]
		\centering
		\begin{minipage}{0.48\textwidth}
			\centering
			\includegraphics[width=\linewidth]{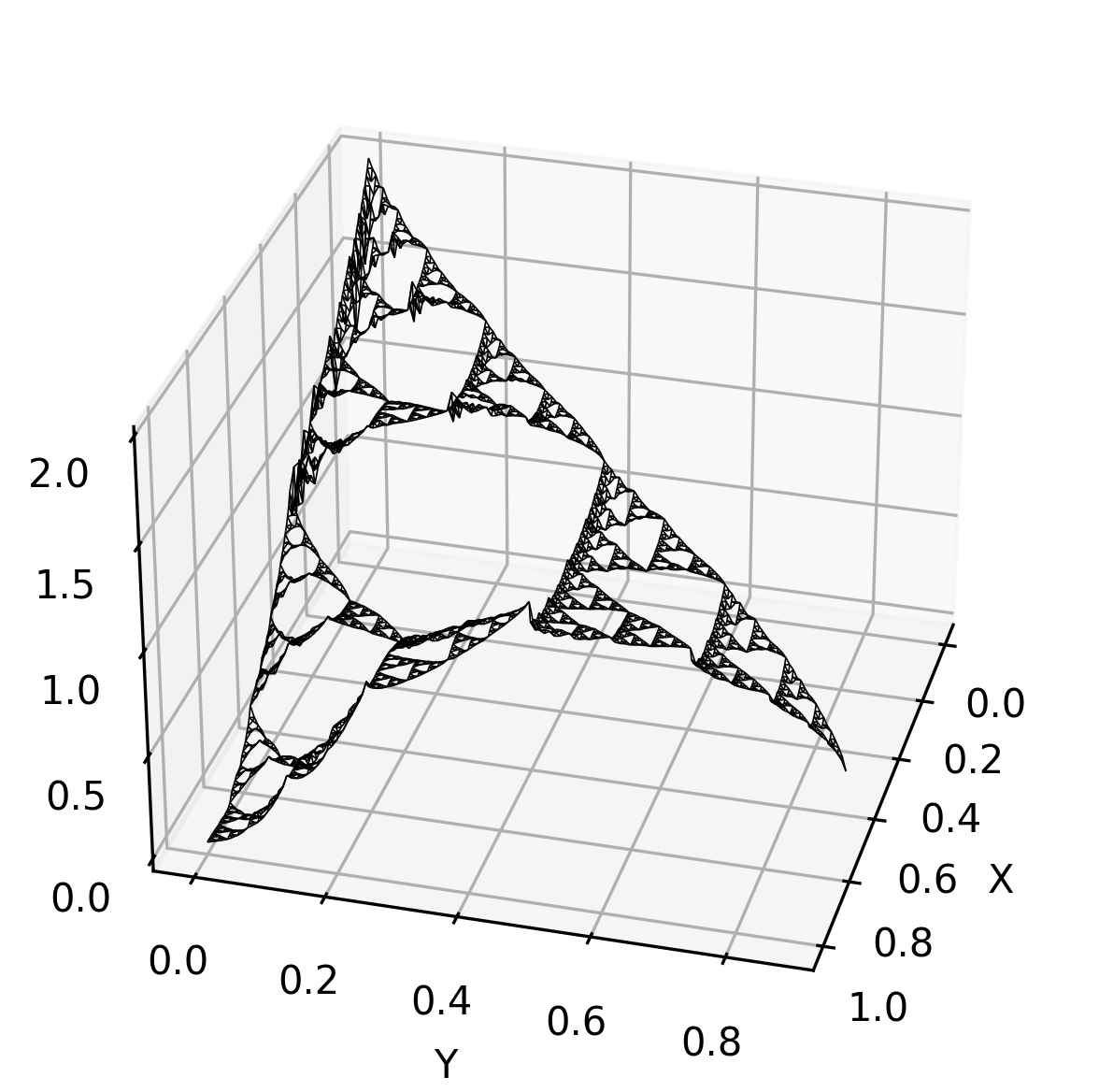}
			\caption{A nonlinear FIF on SG with $\mu<1.5$}
			\label{fig:dim184}
		\end{minipage}
		\hfill
		\begin{minipage}{0.48\textwidth}
			\centering
			\includegraphics[width=\linewidth]{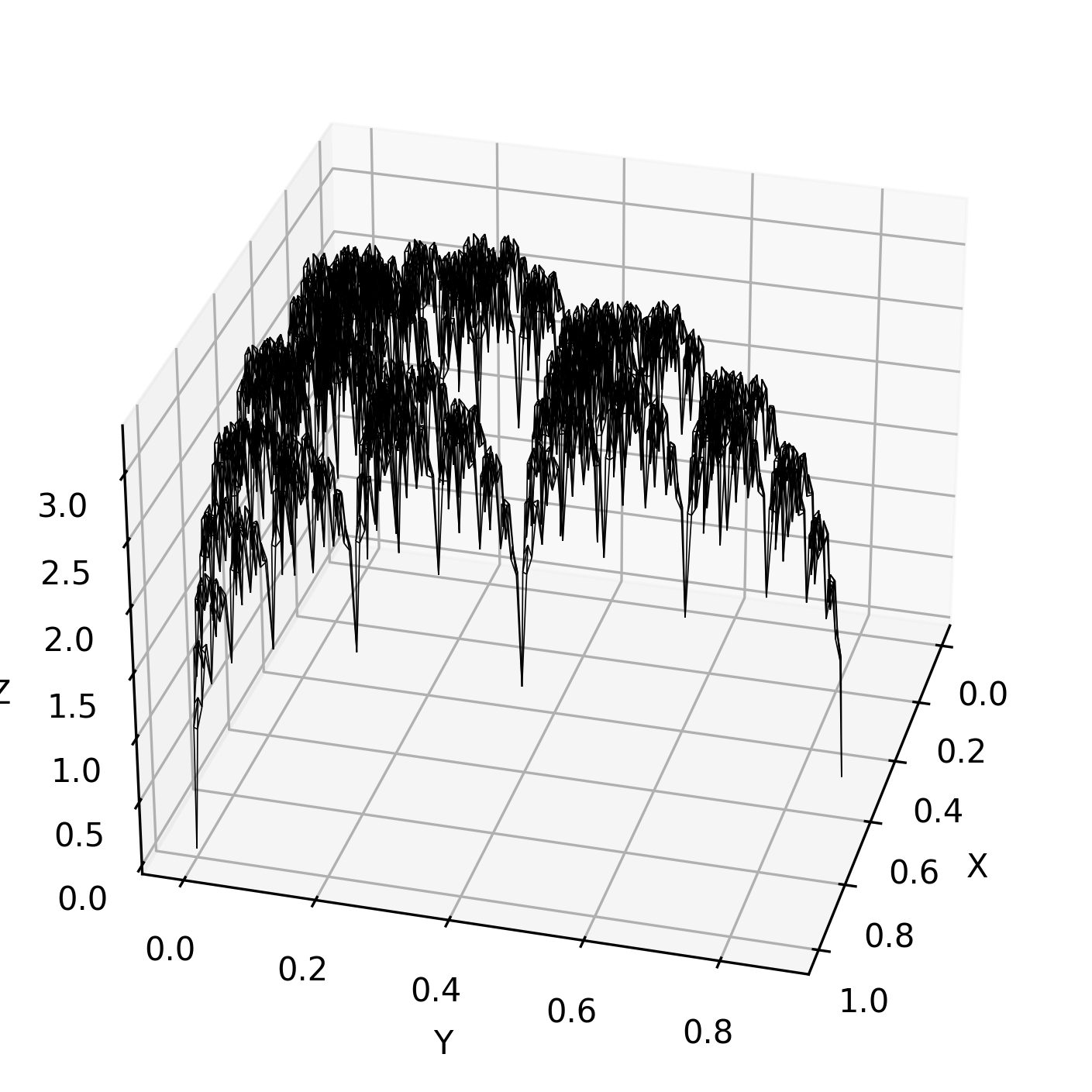}
			\caption{A nonlinear FIF on SG with $\mu=2.7$ and $\nu=2.4$}
			\label{fig:dim224}
		\end{minipage}
		
	\end{figure}

	In the following examples, we construct nonlinear interpolants on the Sierpi\'nski gasket using the prescribed interpolation data and different families of maps $s_n$. We then calculate the corresponding upper and lower bounds for their box-counting dimensions. 
	\begin{example}\label{SGex1}
		For $n\in \{1,2,3\}$, let $l_n:\mathbb{R}^2\to \mathbb{R}^2$ be given as $l_n=\frac{1}{2}(x-q_n)+q_n$, where $q_1=(0,0),q_2=(1,0), q_3=(0.5,\frac{\sqrt{3}}{2})$ and $s_n:[0,2]\to [0,2]$ be defined as
		\begin{align*}
			s_1(y)=\frac{4}{5}\left(\frac{1+y}{2+y}\right),\qquad s_2(y)=\frac{1}{4}y^2,\qquad s_3(y)=\frac{1}{5}y.
		\end{align*}
		Consider the interpolation data as\[
		\begin{aligned}
			&\left\{ (0, 0, 2), \left(0.5, 0, 1\right), (1, 0, 0),\left(0.25, \frac{\sqrt{3}}{4}, 1\right), \left(0.75, \frac{\sqrt{3}}{4}, 1\right), \left(\frac{0.5,\sqrt{3}}{2}, 0\right) \right\}.
		\end{aligned}
		\] The data obeys the following functional equation
		
		$$f^*(l_{n}(x)) = s_{n}(f^*(x)) + h_{n}(x)\\ \,  \text{for all}\;\, \ n \in \{1,2,3\} \;\text{and}\ x \in \mathcal{K},$$
		where $h_n$ is the harmonic function constructed through the given interpolation data and the choice of $s_n$. That is,
		\begin{align*}
			&h_1(q_1)=\frac75,\quad
			h_1(q_2)=\frac35,\quad h_1(q_3)=\frac35,\\
			&h_2(q_1)=0,\quad
			h_2(q_2)=0,\quad h_2(q_3)=1,\\
			&h_3(q_1)=\frac35,\quad
			h_3(q_2)=1,\quad h_3(q_3)=0.
		\end{align*}
		
		Now, $\mu=L(s_1)+L(s_2)+L(s_3)=\frac{1}{5}+1+\frac{1}{5}=1.4<3/2.$
		Hence, from Theorem \ref{SGUpperBox}
		$ \dim_B(G)=\frac{\log3}{\log2}\approx1.585. 
		$
		 Figure \ref{fig:dim184} represents the corresponding FIF.
	\end{example}
	\begin{example}\label{SGex2}
		Consider the interpolation data
		\[
		u(q_i)=0,\qquad
		u\left(\frac{q_i+q_j}{2}\right)=1
		\quad\text{for }i\ne j,
		\]
		and let $s_n:[0,10]\rightarrow[0,10]$ be defined by
		\[
		s_n(y)=\frac45y+\frac1{200}y^2,
		\qquad n=1,2,3.
		\]
		Let $h_n$ be the harmonic functions constructed
		through the given interpolation data and the
		choice of $s_n$. That is,
		\begin{align*}
			&h_1(q_1)=0,\quad h_1(q_2)=1,\quad h_1(q_3)=1,\\
			&h_2(q_1)=1,\quad h_2(q_2)=0,\quad h_2(q_3)=1,\\
			&h_3(q_1)=1,\quad h_3(q_2)=1,\quad h_3(q_3)=0.
		\end{align*}
		Since $0\leq h_n\leq1$ and
		\[
		0\leq h_n(x)+s_n(y)
		\leq1+s_n(10)=\frac{19}{2}<10,
		\]
		the full vertical maps preserve $[0,10]$. On this interval,
		\[
		s_n'(y)=\frac45+\frac{y}{100},
		\]
		so $l(s_n)=4/5$ and $L(s_n)=9/10$.
		Therefore,
		\begin{align*}
			\nu&=l(s_1)+l(s_2)+l(s_3)
			=\frac{12}{5}>\frac32,\\
			\mu&=L(s_1)+L(s_2)+L(s_3)
			=\frac{27}{10}>\frac32.
		\end{align*}
		
		Furthermore, the harmonic extension rule gives
		\[
		H_1=3\left(\frac35+\frac25+\frac25\right)
		=\frac{21}{5}.
		\]
		For each $n$, let $r_n$ be the midpoint of the edge
		opposite $q_n$. Since $f^*(r_n)=1$ and
		$h_n(r_n)=4/5$, we obtain
		\[
		f^*(l_n(r_n))
		=s_n(1)+\frac45
		=\frac{321}{200}.
		\]
		Since $f^*(q_n)=0$ and $q_n\in l_n(\mathcal K)$,
		\[
		S_1\geq\frac{963}{200}>
		\frac{21}{5}
		=\frac{H_1}{\nu-\frac75}.
		\]
		Hence, Theorems \ref{SGUpperBox} and \ref{SGUpperBox} give
		\[
		2.263\approx1+\frac{\log(12/5)}{\log2}
		\leq\underline{\dim}_B(G)
		\leq\overline{\dim}_B(G)
		\leq1+\frac{\log(27/10)}{\log2}\approx2.433.
		\]
		Figure \ref{fig:dim224} represents the corresponding FIF.
	\end{example}
	%
	\subsection{Box dimensions of nonlinear FIFs on KC}
	The results on box dimensions of nonlinear FIFs on the von Koch curve require lemmas similar to the ones shown in the case of the Sierpi\'nski gasket.
	

	\begin{lemma}\label{KCcover}
	Let $\mathcal{V}$ be the standard von Koch curve
	with endpoints $(0,0)$ and $(1,0)$, and let
	$f^*:\mathcal{V}\to\mathbb{R}$ be continuous.
	Write
	\[
	G=\{(x,f^*(x)):x\in\mathcal{V}\},
	\qquad
	G_w=\{(x,f^*(x)):x\in l_w(\mathcal{V})\}.
	\]
	Let $N(k)$ and $N_w(k)$ denote the minimum numbers
	of cubes of side length $3^{-k}$ required to
	cover $G$ and $G_w$, respectively. Then
	\[
	\frac1{63}\sum_{|w|=k}N_w(k)
	\leq N(k)
	\leq\sum_{|w|=k}N_w(k).
	\]
	\end{lemma}
	
	\begin{proof}
	Combining minimal covers of the cell graphs gives
	\[
	N(k)\leq\sum_{|w|=k}N_w(k).
	\]
	
	For the reverse inequality, put $\varepsilon=3^{-k}$
	and choose a minimal cover $\{Q_j\}_{j=1}^{N(k)}$
	of $G$ by cubes of side length $\varepsilon$.
	As in the proof of
	Lemma~\ref{SGcover},
	\[
	\sum_{|w|=k}N_w(k)
	\leq\sum_{j=1}^{N(k)}
	\#\{w:|w|=k,\ Q_j\cap G_w\neq\varnothing\}.
	\]
	
	Let $T$ be the closed triangle with vertices
	\[
	(0,0),\qquad (1,0),\qquad
	\left(\frac12,\frac{\sqrt3}{6}\right).
	\]
	For the standard Koch maps, $\mathcal{V}\subset T$,
	the triangles $l_i(T)$ lie in $T$, and their
	interiors are pairwise disjoint. Consequently,
	the triangles $T_w=l_w(T)$, $|w|=k$, also have
	pairwise disjoint interiors. Each has a diameter
	$\varepsilon$ and an area
	\[
	|T_w|=\frac{\sqrt3}{12}\varepsilon^2.
	\]
	
	If $Q_j$ meets $G_w$, its horizontal projection
	meets $T_w$. Every such triangle is contained
	in the square obtained by extending that
	projection by $\varepsilon$ in each coordinate
	direction. This square has an area
	of $9\varepsilon^2$. Hence, the number of such
	triangles is at most
	\[
	\frac{9\varepsilon^2}
	{\sqrt3\varepsilon^2/12}
	=36\sqrt3<63.
	\]
	Therefore,
	\[
	\sum_{|w|=k}N_w(k)\leq63N(k).
	\]
	\end{proof}
	\begin{lemma}
		\label{KCharmosc}
		Let  $h$ be a harmonic function on the von Koch curve $\mathcal{V}$.  Then, for every	$w\in W_*$,
		\[
		\sum_{i=1}^{4}
		\operatorname{osc}_{l_{wi}(\mathcal{V})}h
		=
		\operatorname{osc}_{l_w(\mathcal{V})}h.
		\]
	\end{lemma}
	
	\begin{proof}
		Let $z_0,z_1,\ldots,z_4$ be the consecutive first-level
		vertices of $\mathcal{V}$, where $z_{i-1}$ and $z_i$
		are the endpoints of $l_i(\mathcal{V})$. Write
		\[
		h(z_0)=A,\qquad h(z_4)=B.
		\]
		The harmonic extension minimizes
		\[
		\sum_{i=1}^{4}(v_i-v_{i-1})^2
		\]
		subject to $v_0=A$ and $v_4=B$. By the
		Cauchy-Schwartz inequality
		\[
		\sum_{i=1}^{4}(v_i-v_{i-1})^2
		\geq
		\frac14\left(
		\sum_{i=1}^{4}(v_i-v_{i-1})
		\right)^2
		=
		\frac{(B-A)^2}{4}.
		\]
		Equality holds when all four increments
		are equal. Hence,
		\[
		h(z_i)=A+\frac{i}{4}(B-A),
		\qquad i=0,\ldots,4.
		\]
		By the harmonic maximum principle
		\[
		\operatorname{osc}_{l_i(\mathcal{V})}h
		=
		|h(z_i)-h(z_{i-1})|
		=
		\frac14|B-A|.
		\]
		Consequently,
		\[
		\sum_{i=1}^{4}
		\operatorname{osc}_{l_i(\mathcal{V})}h
		=
		|B-A|
		=
		\operatorname{osc}_{\mathcal{V}}h.
		\]
		Finally, applying this identity to the harmonic
		function $h\circ l_w$ proves the assertion
		for every $w\in W_*$.
	\end{proof}
	We obtain the following results by a similar proof done in Theorems \ref{SGUpperBox} and \ref{SGLowerBox}, using Lemmas \ref{KCcover} and \ref{KCharmosc}.
\begin{theorem}
	\label{KCUpperBox}
	Let $G=\{(x,f^*(x)):x\in\mathcal{V}\}$ be the graph
	of $f^*$ defined in Theorem~\ref{graph}. Define
	$\mu=\sum_{n=1}^{4}L(s_n)$.
	Then
	\[
	\overline{\dim}_B(G)\leq
	\begin{cases}
		\frac{\log4}{\log3},
		&\mu\leq\frac43,\\[4pt]
		1+\frac{\log\mu}{\log3},
		&\mu>\frac43.
	\end{cases}
	\]
	In particular, if $\mu\leq4/3$, then
	\[
	\dim_B(G)=\frac{\log4}{\log3}.
	\]
\end{theorem}

\begin{theorem}
	\label{KCLowerBox}
	Let $G=\{(x,f^*(x)):x\in\mathcal{V}\}$ be the graph
	of $f^*$ defined in Theorem~\ref{graph}. Define
	$\nu=\sum_{n=1}^{4}l(s_n)$, and let
	\[
	S_1=\sum_{i=1}^{4}
	\operatorname{osc}_{l_i(\mathcal V)}f^*,
	\qquad
	H_1=\sum_{n=1}^{4}\sum_{i=1}^{4}
	\operatorname{osc}_{l_i(\mathcal V)}h_n.
	\]
	If $\nu>\frac43$ and
	$S_1>\frac{H_1}{\nu-1}$, then
	\[
	\underline{\dim}_B(G)
	\geq1+\frac{\log\nu}{\log3}.
	\]
\end{theorem}

\begin{figure}[h]
		\centering
		\begin{minipage}{0.48\textwidth}
			\centering
			\includegraphics[width=\linewidth]{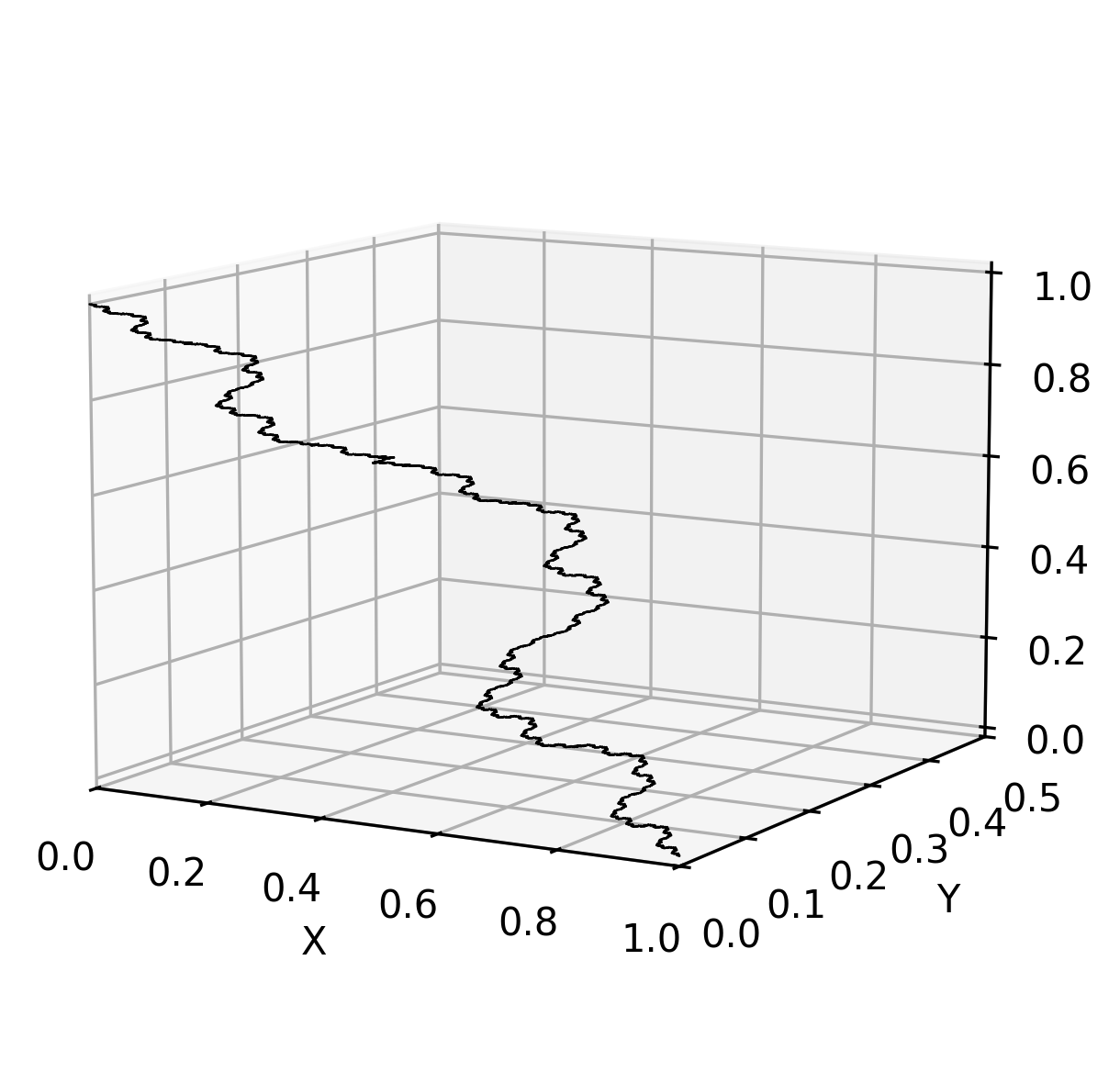}
			\caption{A harmonic function on KC}
			\label{fig:kcharm}
		\end{minipage}
		\hfill
		\begin{minipage}{0.48\textwidth}
			\centering
			\includegraphics[width=\linewidth]{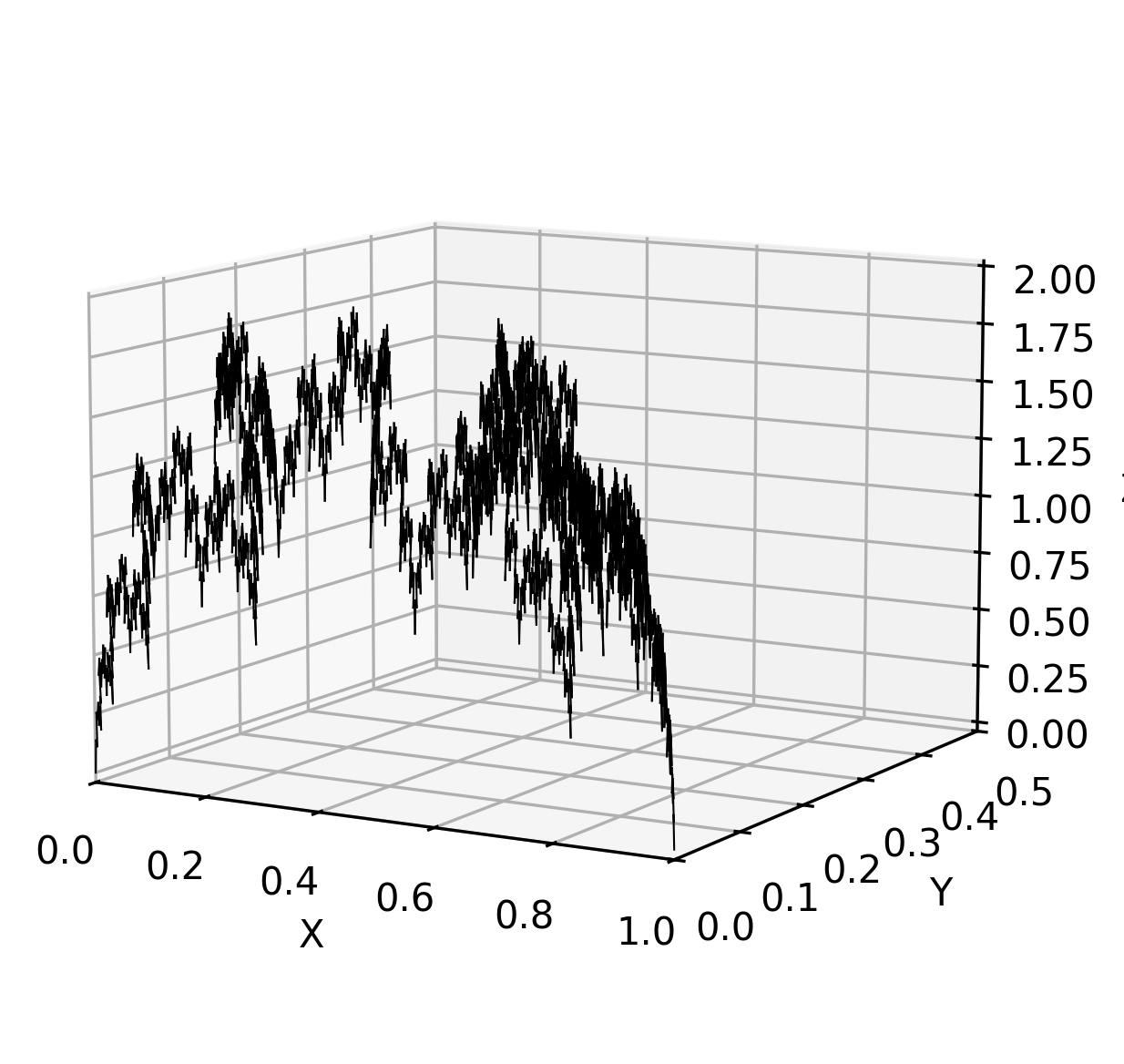}
			\caption{A nonlinear FIF on KC with $\mu=2.64$}
			\label{fig:dim188}
		\end{minipage}
		
	\end{figure}

	\begin{corollary}
		\label{cor:harmonic-dimension-KC}
		Let $h:\mathcal{V}\to\mathbb{R}$ be harmonic on the
		standard von Koch curve equipped with its standard
		symmetric harmonic structure. Then, the graph $
		G_h=\{(x,h(x)):x\in\mathcal{V}\}
		$
		has a box-counting dimension
		$
		\dim_B(G_h)=\frac{\log4}{\log3}.
		$
	\end{corollary}
	\begin{example}\label{Ex-KC-lower}
		Let $z_0,z_1,\ldots,z_4$ be the consecutive
		first-level vertices of the von Koch curve
		$\mathcal V$, where $z_0=q_1$ and $z_4=q_2$.
		Consider the interpolation data
		$
		u(z_0)=u(z_4)=0,\;
		u(z_1)=u(z_2)=u(z_3)=1,
		$
		and let $s_n:[0,3]\rightarrow[0,3]$ be defined by
		\[
		s_n(y)=\frac35y+\frac1{100}y^2,
		\; n=1,2,3,4.
		\]
		Let $h_n$ be the harmonic functions constructed
		through the given interpolation data and the
		choice of $s_n$. That is, $
			h_1(q_1)=0,\; h_1(q_2)=1,
			h_2(q_1)=1,\; h_2(q_2)=1,
			h_3(q_1)=1,\; h_3(q_2)=1,
			h_4(q_1)=1,\; h_4(q_2)=0.
		$
		Since $0\leq h_n\leq1$, and
		$
		0\leq h_n(x)+s_n(y)
		\leq1+s_n(3)=\frac{289}{100}<3,
		$
		the full vertical maps preserve $[0,3]$. On this interval,
		\[
		s_n'(y)=\frac35+\frac{y}{50},
		\]
		so $l(s_n)=3/5$ and $L(s_n)=33/50$.
		Therefore,
		\begin{align*}
			\nu&=\sum_{n=1}^{4}l(s_n)
			=4\cdot\frac35
			=\frac{12}{5}>\frac43,\\
			\mu&=\sum_{n=1}^{4}L(s_n)
			=4\cdot\frac{33}{50}
			=\frac{66}{25}>\frac43.
		\end{align*}
		
		By the harmonic oscillation identity,
		\[
		H_1=\sum_{n=1}^{4}
		\operatorname{osc}_{\mathcal V}h_n
		=1+0+0+1=2.
		\]
		The interpolation values at the endpoints of
		the first and fourth cells give
		\[
		\operatorname{osc}_{l_1(\mathcal V)}f^*\geq1,
		\qquad
		\operatorname{osc}_{l_4(\mathcal V)}f^*\geq1.
		\]
		Consequently,
		\[
		S_1\geq2>\frac{10}{7}
		=\frac{H_1}{\nu-1}.
		\]
		Hence, from Theorems \ref{KCUpperBox} and \ref{KCLowerBox} give
		\[
		1.797\approx1+\frac{\log(12/5)}{\log3}
		\leq\underline{\dim}_B(G)
		\leq\overline{\dim}_B(G)
		\leq1+\frac{\log(66/25)}{\log3}\approx1.884.
		\]
		Figure \ref{fig:dim188} represents the corresponding FIF.
	\end{example}
\begin{remark}
	In 2025, Liang and Ruan \cite{Liang} proved that the box
	and Hausdorff dimensions of the graph of a harmonic function on the $(N-1)$-dimensional Sierpi\'nski gasket are both $\log_2 N$. Our results recover the	box-counting dimension formula when $N=3$. We also establish the corresponding equality for the graph of a harmonic function on the standard von Koch curve, that is,
	$\dim_B(G_h)=\log_3 4$.
\end{remark}
	%
	\section{Properties of uniform nonlinear fractal functions}\label{sec-5}
	In 2011, Ri and Ruan \cite{Ri2011} investigated several properties of uniform fractal interpolation functions with a fixed vertical scaling factor on the Sierpi\'nski gasket. 
	 In this section, we establish several properties of uniform nonlinear fractal interpolation functions on the Sierpi\'nski gasket and the von Koch curve.
	\subsection{Properties on SG}
	We call a nonlinear fractal interpolation function uniform if its interpolation data are $0$ on $V_0$ and $1$ on $V_1\backslash V_0$, and if 
	$s_n(0)=0$ for $n=1,2,3$. In this section, consider $V_0=\{q_1,q_2,q_3\}$ for $\mathcal{K}$ and $q_\omega=l_{\omega_1\omega_2\ldots\omega_{k-1}}(q_{\omega_k})$.
	\begin{definition}(Graph energy on SG)\cite{Strichartz}
		For $m\in\mathbb{N}_0$, we define graph energies
		$$\mathcal{E}_m(u):=\left(\frac{5}{3}\right)^m\sum_{x\sim_my}\left(u(x)-u(y)\right)^2,$$
		where $u:V^*\rightarrow\mathbb{R}$, and $x\sim_my$ implies that $x$ and $y$ are neighbours in $V_m$ that is $x\sim_my$ if there exists
		 $\; \omega\in\Sigma^m$ such that $x,y\in l_\omega(V_0)$.
	\end{definition}    
	The graph energy sequence is increasing, and it satisfies the following relation:
	$$\mathcal{E}_{m-1}(u)=\min\mathcal{E}_m(\tilde u),$$
	where the minimum is taken over all $\tilde u$ which satisfies $\tilde u|_{V_{m-1}}=u$. We define the energy of the function as
	$$\mathcal{E}(u):=\lim_{m\rightarrow\infty} \mathcal{E}_{m}(u).$$
	If the limit exists, we say the function $u$ has finite energy.
	\begin{theorem}\label{Thm-energySC}
		Let $f^*$ be the uniform nonlinear FIF on SG and let $\Psi=\sum_{n=1}^3(L(s_n))^2$.  If $\Psi<\frac{3}{5}$, then $\mathcal{E}(f^*)<\infty$. In particular, $\mathcal{E}(f^*)<\frac{30}{3-5\Psi}$.
	\end{theorem}
	\begin{proof}
		From the definition of graph energies for $k\geq2$,
		\allowdisplaybreaks
		\begin{align}\label{eq1-new}
			\mathcal{E}_k(f^*)&=\left(\frac{5}{3}\right)^k\sum_{x\sim_ky}\left(f^*(x)-f^*(y)\right)^2 \nonumber \\
			&=\left(\frac{5}{3}\right)^k\sum_{\omega\in\Sigma^{k-1}}\sum_{1\leq i<j\leq3}\sum_{n=1}^3\left(f^*(q_{n\omega i})-f^*(q_{n\omega j} )\right)^2.
		\end{align}
		For any $\omega\in \Sigma^{k-1}$, and for $i\neq j$, and using Theorem \ref{graph}, we get
		\allowdisplaybreaks
		\begin{align*}
			&\sum_{n=1}^3\left(f^*(q_{n\omega i})-f^*(q_{n\omega j} )\right)^2= \sum_{n=1}^3\left(s_n(f^*(q_{\omega i}))+h_n(q_{\omega i} )-(s_n(f^*(q_{\omega j}))+h_n(q_{\omega j} ))\right)^2\\
			&=\sum_{n=1}^3\left[(s_n(f^*(q_{\omega i}))-s_n(f^*(q_{\omega j}))]+[h_n(q_{\omega i} )-h_n(q_{\omega j} )\right)]^2\\
			&=\sum_{n=1}^3(s_n(f^*(q_{\omega i}))-s_n(f^*(q_{\omega j}))^2+\sum_{n=1}^3 (h_n(q_{\omega i} )-h_n(q_{\omega j} ))^2\\
			&+ 2\sum_{n=1}^3(s_n(f^*(q_{\omega i}))-s_n(f^*(q_{\omega j}))\cdot(h_n(q_{\omega i} )-h_n(q_{\omega j} ))\\
			&\leq\sum_{n=1}^3(L(s_n))^2(\left(f^*(q_{\omega i})-f^*(q_{\omega j} )\right)^2+\sum_{n=1}^3 (h_n(q_{\omega i} )-h_n(q_{\omega j} ))^2\\
			&+2 \max_{1\leq n\leq 3}(L(s_n))((f^*(q_{\omega i}))-(f^*(q_{\omega j}))\cdot\sum_{n=1}^3(h_n(q_{\omega i} )-h_n(q_{\omega j} ))\\
			&=\Psi\left(f^*(q_{\omega i})-f^*(q_{\omega j} )\right)^2+\sum_{n=1}^3 (h_n(q_{\omega i} )-h_n(q_{\omega j} ))^2\\
			&+2 \max_{1\leq n\leq 3}(L(s_n))((f^*(q_{\omega i}))-(f^*(q_{\omega j}))\cdot\sum_{n=1}^3(h_n(q_{\omega i} )-h_n(q_{\omega j} )).
		\end{align*}
		Since $s_n(0)=0$, we have that $h_n(q_i)=1-\delta_{ni}$ for $n,i=1,2,3$. Hence, $\sum_{n=1}^3h_n(q_i)=2$. Therefore, $\sum_{n=1}^3h_n$ take constant values for points $l_\omega(q_i)$. This implies that $\sum_{n=1}^3h_n(q_{\omega i})-h_n(q_{\omega j})=0$, and hence the last term in the previous inequality vanishes. From Equation \eqref{eq1-new}, we obtain
		\begin{align}\label{inequa}
			\mathcal{E}_k(f^*)&\leq\left(\frac{5}{3}\right)^k\Psi\sum_{\omega\in\Sigma^{k-1}}\sum_{1\leq i<j\leq3}\left(f^*(q_{\omega i})-f^*(q_{\omega j} )\right)^2\nonumber\\
			&+\left(\frac{5}{3}\right)^k\sum_{\omega\in\Sigma^{k-1}}\sum_{1\leq i<j\leq3}\sum_{n=1}^3\left(h_n(q_{\omega i})-h_n(q_{\omega j} )\right)^2\nonumber\\
			&=\frac{5\Psi}{3}\mathcal{E}_{k-1}(f^*)+\frac{5}{3}\sum_{n=1}^3\mathcal{E}_{k-1}(h_n).
		\end{align}
		Since $h_n$ are harmonic, $\mathcal{E}_{k-1}(h_n)=\mathcal{E}_0(h_n)=1^2+1^2+0^2=2$. Since $f^*$ is defined on $V_1$, $\mathcal{E}_1(f^*)=\frac{5}{3}(2+2+2)=10$. Substituting these values in Inequality \eqref{inequa}, we get
		\begin{align*}
			\mathcal{E}_k(f^*)&\leq\frac{5\Psi}{3}\mathcal{E}_{k-1}(f^*)+10\\
			&\leq\frac{5\Psi}{3}\left(\frac{5\Psi}{3}\mathcal{E}_{k-2}(f^*)+10\right)+10\\
			&\vdots\\
			&\leq \left(\frac{5\Psi}{3}\right)^{k-1}\mathcal{E}_1(f^*)+10\left(1+\frac{5\Psi}{3}+\ldots+\left(\frac{5\Psi}{3}\right)^{k-2}\right)\\
			&=10\left(\frac{5\Psi}{3}\right)^{k-1}+ 10\left(\frac{(\frac{5\Psi}{3})^{k-1}-1}{\frac{5\Psi}{3}-1}\right).
		\end{align*}
		Assume $\Psi<\frac{3}{5}$. Then
		\begin{align*}
			\mathcal{E}(f^*)&=\lim_{k\rightarrow\infty} \mathcal{E}_{k}(u)\\
			&\leq\lim_{k\rightarrow\infty} \left(10\left(\frac{5\Psi}{3}\right)^{k-1}+ 10\left(\frac{(\frac{5\Psi}{3})^{k-1}-1}{\frac{5\Psi}{3}-1}\right)\right)\\
			&=\lim_{k\rightarrow\infty} 10\left(\frac{5\Psi}{3}\right)^{k-1}+ \lim_{k\rightarrow\infty}10\left(\frac{(\frac{5\Psi}{3})^{k-1}-1}{\frac{5\Psi}{3}-1}\right)\\
			&=\frac{30}{3-5\Psi}.
		\end{align*}
		Hence the proof.
		\begin{figure}[h]
			\centering
			\includegraphics[width=0.5\linewidth]{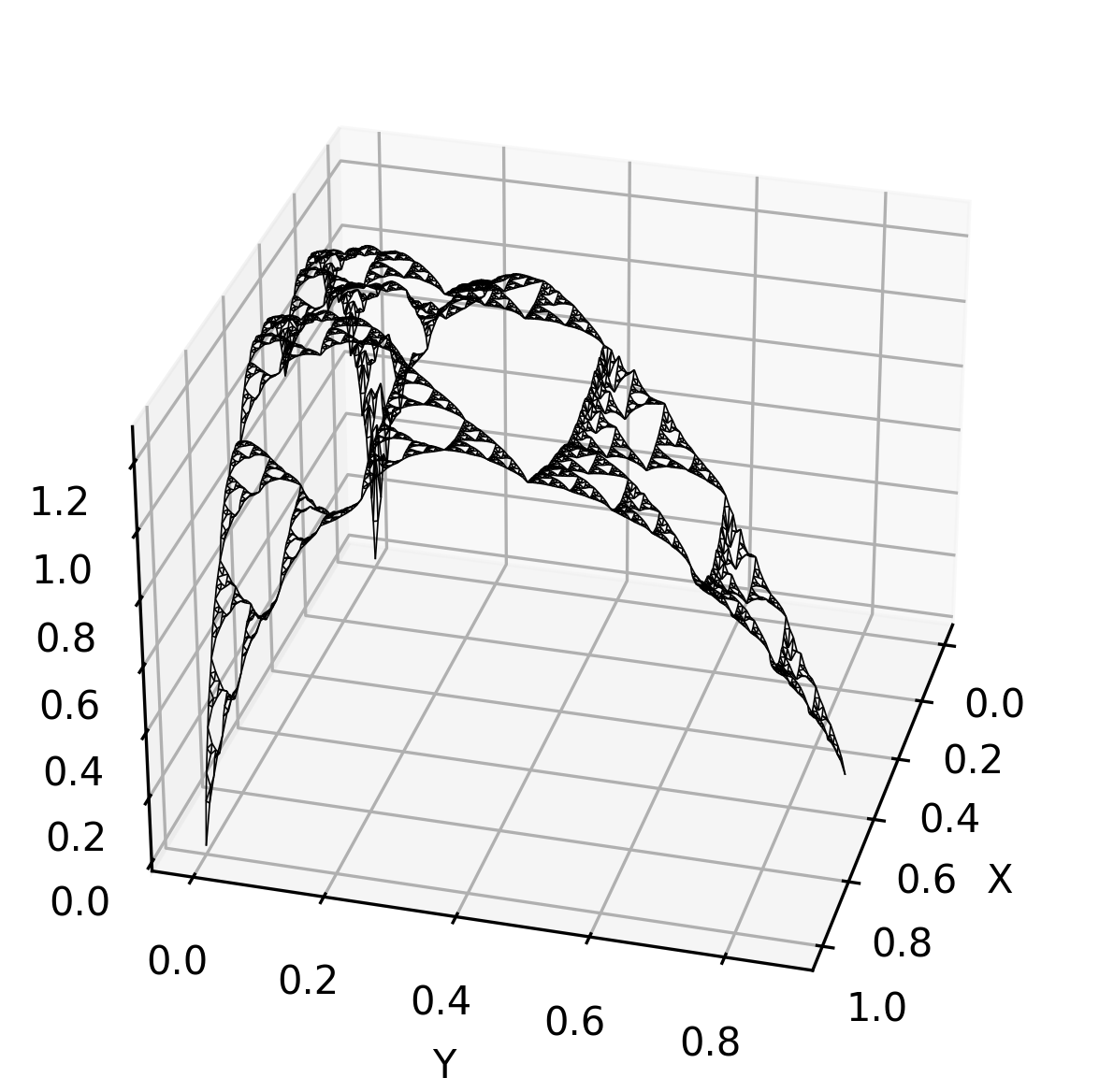}
			\caption{A nonlinear FIF on SG with finite energy ($\Psi=0.515$)}
			\label{fig:finenerg}
		\end{figure}
	\end{proof}
	\begin{definition}(Normal derivative on SG)
		\cite{Strichartz}	The normal derivative of $u$ at a point $q_i\in V_0$ is defined as
		$$\partial_n(u(q_i)):=\lim_{m\rightarrow\infty}\left(\frac{5}{3}\right)^m(2u(q_i)-u(F_i^m(q_j))-u(F_i^m(q_k)))$$
		when the limit exists, where $i,j,k$ are distinct.
	\end{definition}
	It is easy to see that for harmonic functions, the normal derivative is always well defined and is given by $\partial_n(h(q_i))=2h(q_i)-h(q_j)+h(q_k)$.
	Now we give a theorem for the existence of the normal derivative of the uniform nonlinear FIF.
	\begin{lemma}
		Let $f^*$ be the uniform nonlinear FIF such that $s_n$ is an increasing function for $n=1,2,3$. Then $f^*(x)\geq0$ for all $x\in \mathcal{K}$.
	\end{lemma}
	\begin{proof}
		We prove this by induction. Since we know that $s_n$ is increasing for $n=1,2,3$, we have
		$$s_i(x)-s_i(0)\geq l(s_i)\cdot x ,$$
		for any $x>0.$ Clearly, $h_n$ is a non-negative function. Now, assume that $f^*(q_\omega)\geq0$ for all $\omega\in \Sigma^k$. Then by Theorem \ref{graph},
		\begin{align*}
			f^*(q_{n\omega})&=s_n(f^*(q_\omega))+h_n(q_\omega)\\
			&\geq l(s_n) f^*(q_\omega)+h_n(q_\omega).
		\end{align*}
		By our assumption, $f^*(q_\omega)\geq0$ and $h_n(q_\omega)\geq0$. Then $f^*(q_\omega)\geq 0$.
		Also, we have $f^*(q_\omega)\geq0$ for $\omega\in \Sigma^2$. Therefore, by induction, $f^*(q_\omega)\geq0$ for $\omega\in \Sigma^k$. Using the continuity of $f^*$, we get $f^*(x)\geq0$ for all $x\in \mathcal{K}$.
	\end{proof}
	\begin{theorem}\label{Thm-Normal}
		Let $f^*$ be the non-negative uniform nonlinear FIF such that $s_i$ is an increasing function for $i=1,2,3$. Then for $i\in \{1,2,3\}$, the normal derivative $\partial_n(q_i)$ exists if $L(s_i)<3/5$. Furthermore, the bounds on the normal derivative are given by
		$$\frac{2}{L(s_i)-3/5}\leq\partial_nf^*(q_i)\leq\frac{2}{l(s_i)-3/5}.$$
	\end{theorem}
	\begin{proof}
		Let $i\neq j$. We denote $[i]^m=\underbrace{ii\ldots i}_{m\ \text{times}}$.
		Applying the $\frac{2}{5}$-$\frac{1}{5}$ rule to the harmonic function $h_i,$ we obtain $h_i(q_{[i]^mj})=\left(\frac{3}{5}\right)^m$.
		By Theorem \ref{graph}, we have
		$$f^*(q_{[i]^{m+1}j})=s_i(f^*(q_{[i]^mj})+h_i(q_{[i]j})=s_i(f^*(q_{[i]^mj})+\left(\frac{3}{5}\right)^m.$$
		Since $s_i(0)=0,\; |s_i(x)-s_i(0)|\leq L(s_i)|x-0| \implies |s_i(x)|\leq L(s_i)|x|\; \text{for all} \; x\in\mathcal{K}$. Hence, we have
		\begin{align*}
			f^*(q_{[i]^{m+1}j})& \leq L(s_i)(f^*(q_{[i]^mj})+\left(\frac{3}{5}\right)^m.
		\end{align*}
		Using this inequality repeatedly and the fact that $f^*(q_j)=0$, we get
		\begin{align*}
			f^*(q_{[i]^{m+1}j}) & \leq L(s_i)(f^*(q_{[i]^mj})+\left(\frac{3}{5}\right)^m\\
			&\leq L(s_i)\left(L(s_i)f^*(q_{[i]^{m-1}j}) +\left(\frac{3}{5}\right)^{m-1}\right) +\left(\frac{3}{5}\right)^m\\
			&\vdots\\
			& \leq \sum_{k=0}^m(L(s_i))^k \left(\frac{3}{5}\right)^{m-k}  \\
			&= \frac{L(s_i)^{m+1}-(3/5)^{m+1}}{L(s_i)-3/5}
		\end{align*}
		when $L(s_i)\neq3/5$.
		Similarly, we get 
		$$f^*(q_{[i]^{m+1}k}) \leq \frac{L(s_i)^{m+1}-(3/5)^{m+1}}{L(s_i)-3/5},$$
		for $k\neq i$. Assume $L(s_i)<3/5$. 
		Now, by using the definition of normal derivative, we get
		\begin{align*}
			\partial_n (f^*(q_i))&=\lim_{m\rightarrow\infty}\left(\frac{5}{3}\right)^m(2f^*(q_i)-f^*(F_i^m(q_j))-f^*(F_i^m(q_k)))\\
			&=\lim_{m\rightarrow\infty}\left(\frac{5}{3}\right)^m(2f^*(q_i)-f^*(q_{[i]^mj})-f^*(q_{[i]^mk}))\\
			&=-\lim_{m\rightarrow\infty}\left(\frac{5}{3}\right)^m f^*(q_{[i]^mj}) - \lim_{m\rightarrow\infty}\left(\frac{5}{3}\right)^m f^*(q_{[i]^mk})\\
			&  \geq-2\lim_{m\rightarrow\infty}\left(\frac{5}{3}\right)^m \frac{L(s_i)^{m}-(3/5)^{m}}{L(s_i)-3/5}             \\
			&\geq \frac{2}{L(s_i)-3/5}.
		\end{align*}
		A similar argument yields
		$$f^*(q_{[i]^{m+1}k}) \geq \frac{l(s_i)^{m+1}-(3/5)^{m+1}}{l(s_i)-3/5} $$
		when $l(s_i)\neq3/5$, and hence,
		\begin{align*}
			\partial_n (f^*(q_i))&=\lim_{m\rightarrow\infty}\left(\frac{5}{3}\right)^m(2f^*(q_i)-f^*(F_i^m(q_j))-f^*(F_i^m(q_k)))\\
			& \leq-2\lim_{m\rightarrow\infty}\left(\frac{5}{3}\right)^m \frac{l(s_i)^{m}-(3/5)^{m}}{l(s_i)-3/5}             \\
			&\leq \frac{2}{l(s_i)-3/5}.
		\end{align*}
	\end{proof}
	\subsection{Properties on KC}
		\begin{definition}(Graph energy on KC)\cite{Frieberg}
		For $m\in\mathbb{N}\cup\{0\}$ we define graph energies $\mathcal{E}_m$ by
		$$\mathcal{E}_m(u):=4^m\sum_{x\sim_my}\left(u(x)-u(y)\right)^2,$$
		where $u:V^*\rightarrow\mathbb{R}$, and $x\sim_my$ implies that $x$ and $y$ are neighbours in $V_m$ that is $x\sim_my$ if there exists $ \omega\in\Sigma^m$ such that $x,y\in l_\omega(V_0)$.
	\end{definition}    
	A similar result on graph energy follows for uniform nonlinear FIFs on the von Koch curve. 
	\begin{theorem}\label{Thm-energyKC}
		Let $f^*$ be the uniform nonlinear FIF on KC and let $\Psi=\sum_{n=1}^4(L(s_n))^2$.  If $\Psi<\frac{1}{4}$, then $\mathcal{E}(f^*)<\infty$. In particular, $\mathcal{E}(f^*)<\frac{24}{1-4\Psi}$.
	\end{theorem}

	\begin{proof}
		From the definition of graph energies for $k\geq2$,
		\allowdisplaybreaks
		\begin{align}\label{eq-new}
			\mathcal{E}_k(f^*)&=4^k\sum_{x\sim_ky}\left(f^*(x)-f^*(y)\right)^2 \nonumber \\
			&=4^k\sum_{\omega\in\Sigma^{k-1}}\sum_{1\leq i<j\leq2}\sum_{n=1}^4\left(f^*(q_{n\omega i})-f^*(q_{n\omega j} )\right)^2.
		\end{align}
		For any $\omega\in \Sigma^{k-1}$, and for $i\neq j$, and using Theorem \ref{graph}, we get
		\allowdisplaybreaks
		\begin{align*}
			&\sum_{n=1}^4\left(f^*(q_{n\omega i})-f^*(q_{n\omega j} )\right)^2= \sum_{n=1}^4\left(s_n(f^*(q_{\omega i}))+h_n(q_{\omega i} )-(s_n(f^*(q_{\omega j}))+h_n(q_{\omega j} ))\right)^2\\
			&=\sum_{n=1}^4\left[(s_n(f^*(q_{\omega i}))-s_n(f^*(q_{\omega j}))]+[h_n(q_{\omega i} )-h_n(q_{\omega j} )\right)]^2\\
			&=\sum_{n=1}^4(s_n(f^*(q_{\omega i}))-s_n(f^*(q_{\omega j}))^2+\sum_{n=1}^4 (h_n(q_{\omega i} )-h_n(q_{\omega j} ))^2\\
			&+ 2\sum_{n=1}^4(s_n(f^*(q_{\omega i}))-s_n(f^*(q_{\omega j}))\cdot(h_n(q_{\omega i} )-h_n(q_{\omega j} ))\\
			&\leq\sum_{n=1}^4(L(s_n))^2(\left(f^*(q_{\omega i})-f^*(q_{\omega j} )\right)^2+\sum_{n=1}^4 (h_n(q_{\omega i} )-h_n(q_{\omega j} ))^2\\
			&+2 \max_{1\leq n\leq 4}(L(s_n))((f^*(q_{\omega i}))-(f^*(q_{\omega j}))\cdot\sum_{n=1}^4(h_n(q_{\omega i} )-h_n(q_{\omega j} ))\\
			&=\Psi\left(f^*(q_{\omega i})-f^*(q_{\omega j} )\right)^2+\sum_{n=1}^4 (h_n(q_{\omega i} )-h_n(q_{\omega j} ))^2\\
			&+2 \max_{1\leq n\leq 4}(L(s_n))((f^*(q_{\omega i}))-(f^*(q_{\omega j}))\cdot\sum_{n=1}^4(h_n(q_{\omega i} )-h_n(q_{\omega j} )).
		\end{align*}
		Since $s_n(0)=0$, we have
		\begin{align*}
			&h_1(q_1)=0, h_1(q_2)=1
			h_2(q_1)=0,h_2(q_2)=1\\
			&h_3(q_1)=0, h_3(q_2)=1
			h_4(q_1)=0, h_4(q_2)=1.
		\end{align*}
		 Hence, $\sum_{n=1}^4h_n(q_i)=3$. This implies that $\sum_{n=1}^4h_n(q_{\omega i})-h_n(q_{\omega j})=0$, and hence the last term in the previous inequality vanishes. From Equation \eqref{eq-new}, we obtain
		\begin{align}\label{inequality}
			\mathcal{E}_k(f^*)&\leq4^k\Psi\sum_{\omega\in\Sigma^{k-1}}\sum_{1\leq i<j\leq2}\left(f^*(q_{\omega i})-f^*(q_{\omega j} )\right)^2\nonumber \\
			&+4^k\sum_{\omega\in\Sigma^{k-1}}\sum_{1\leq i<j\leq2}\sum_{n=1}^4\left(h_n(q_{\omega i})-h_n(q_{\omega j} )\right)^2\nonumber\\
			&=4\Psi\mathcal{E}_{k-1}(f^*)+4\sum_{n=1}^4\mathcal{E}_{k-1}(h_n).
		\end{align}
		Since $h_n$ are harmonic, $\mathcal{E}_{k-1}(h_n)=\mathcal{E}_0(h_n)$. Thus we have
		\begin{align*}
			\mathcal{E}_{k-1}(h_1)=0^2+1^2=1, \mathcal{E}_{k-1}(h_2)=1^2+1^2=2, 		\mathcal{E}_{k-1}(h_3)=1^2+1^2=2, \mathcal{E}_{k-1}(h_4)=1^2+0^2=1.
		\end{align*}
		Since $f^*$ is defined on $V_1$, $\mathcal{E}_1(f^*)=4(1+0+0+1)=8$. Therefore, Inequality \ref{inequality} becomes
		\begin{align*}
			\mathcal{E}_k(f^*)&\leq4\Psi\mathcal{E}_{k-1}(f^*)+24\\
			&\leq4\Psi\left(4\Psi\mathcal{E}_{k-2}(f^*)+24\right)+24\\
			&\vdots\\
			&\leq (4\Psi)^{k-1}\mathcal{E}_1(f^*)+24\left(1+4\Psi+\ldots+\left(4\Psi\right)^{k-2}\right)\\
			&=8(4\Psi)^{k-1}+ 24\left(\frac{(4\Psi)^{k-1}-1}{4\Psi-1}\right).
		\end{align*}
		Assume $\Psi<\frac{1}{4}$. Then
		\begin{align*}
			\mathcal{E}(f^*)&=\lim_{k\rightarrow\infty} \mathcal{E}_{k}(u)\\
			&\leq\lim_{k\rightarrow\infty} \left(8(4\Psi)^{k-1}+ 24\left(\frac{(4\Psi)^{k-1}-1}{4\Psi-1}\right)\right)\\
			&=\lim_{k\rightarrow\infty} 8\left(4\Psi\right)^{k-1}+ \lim_{k\rightarrow\infty}24\left(\frac{(4\Psi)^{k-1}-1}{4\Psi-1}\right)\\
			&=\frac{24}{1-4\Psi}.
		\end{align*}
		Hence the proof.
	\end{proof}
	\section*{Acknowledgments}
	This work is supported by CSIR-HRDG ASPIRE scheme (grant no-25WS(014)/2023-24/EMR-II/ASPIRE). The first and second authors thank CSIR for its financial support.
	\bibliographystyle{amsplain}
	\bibliography{references}

@article {Amo,
	AUTHOR = {de Amo, E. and D\'iaz Carrillo, M. and Fern\'andez
	S\'anchez, J.},
	TITLE = {P{CF} self-similar sets and fractal interpolation},
	JOURNAL = {Math. Comput. Simulation},
	FJOURNAL = {Mathematics and Computers in Simulation},
	VOLUME = {92},
	YEAR = {2013},
	PAGES = {28--39},
}

@book {Barnsley1993,
	AUTHOR = {Barnsley, M. F.},
	TITLE = {Fractals {E}verywhere},
	EDITION = {Second},
	PUBLISHER = {Academic Press Professional, Boston, MA},
	YEAR = {1993},
	PAGES = {xiv+534},
}

@article {Barnsley1986,
	AUTHOR = {Barnsley, M. F.},
	TITLE = {Fractal functions and interpolation},
	JOURNAL = {Constr. Approx.},
	FJOURNAL = {Constructive Approximation. An International Journal for
	Approximations and Expansions},
	VOLUME = {2},
	YEAR = {1986},
	NUMBER = {4},
	PAGES = {303--329},
}

@article {BarnsleyEltonHardin,
	AUTHOR = {Barnsley, M. F. and Elton, J. H. and Hardin, D.
	P.},
	TITLE = {Recurrent iterated function systems},
	JOURNAL = {Constr. Approx.},
	FJOURNAL = {Constructive Approximation. An International Journal for
	Approximations and Expansions},
	VOLUME = {5},
	YEAR = {1989},
	NUMBER = {1},
	PAGES = {3--31},}

@article {Manuj,
    AUTHOR = {Verma, M. and Priyadarshi, A. and Verma, S.},
     TITLE = {Analytical and dimensional properties of fractal interpolation
              functions on the {S}ierpi\'nski gasket},
   JOURNAL = {Fract. Calc. Appl. Anal.},
  FJOURNAL = {Fractional Calculus and Applied Analysis. An International
              Journal for Theory and Applications},
    VOLUME = {26},
      YEAR = {2023},
    NUMBER = {3},
     PAGES = {1294--1325},}

@article {BEH,
	AUTHOR = {Barnsley, M. F. and Elton, J. H. and Hardin, D. and Massopust,
	P.},
	TITLE = {Hidden variable fractal interpolation functions},
	JOURNAL = {SIAM J. Math. Anal.},
	FJOURNAL = {SIAM Journal on Mathematical Analysis},
	VOLUME = {20},
	YEAR = {1989},
	NUMBER = {5},
	PAGES = {1218--1242},
}

@article {BouDalla,
	AUTHOR = {Bouboulis, P. and Dalla, L. and Drakopoulos, V.},
	TITLE = {Construction of recurrent bivariate fractal interpolation
	surfaces and computation of their box-counting dimension},
	JOURNAL = {J. Approx. Theory},
	FJOURNAL = {Journal of Approximation Theory},
	VOLUME = {141},
	YEAR = {2006},
	NUMBER = {2},
	PAGES = {99--117},
}

@article {Celik,
	AUTHOR = {D. \c{C}elik and \c{S}. Ko\c{c}ak and Y. \"Ozdemir},
	TITLE = {Fractal interpolation on the {S}ierpinski gasket},
	JOURNAL = {J. Math. Anal. Appl.},
	FJOURNAL = {Journal of Mathematical Analysis and Applications},
	VOLUME = {337},
	YEAR = {2008},
	NUMBER = {1},
	PAGES = {343--347},
}

@article {Dumitru2009,
	AUTHOR = {Dumitru, D.},
	TITLE = {Generalized iterated function systems containing
	{M}eir-{K}eeler functions},
	JOURNAL = {An. Univ. Bucure\c sti Mat.},
	FJOURNAL = {Analele Universit\u a\c tii Bucure\c sti. Matematic\u a},
	VOLUME = {58},
	YEAR = {2009},
	NUMBER = {1},
	PAGES = {109--121},
}

@article {Edelstein,
	AUTHOR = {Edelstein, M.},
	TITLE = {On fixed and periodic points under contractive mappings},
	JOURNAL = {J. London Math. Soc.},
	FJOURNAL = {The Journal of the London Mathematical Society},
	VOLUME = {37},
	YEAR = {1962},
	PAGES = {74--79},
}

@article{Hu,
	TITLE={Domains of {D}irichlet forms and effective resistance estimates on p.c.f. fractals},
	AUTHOR={Hu, J. and Wang, X.},
	JOURNAL={Studia Mathematica},
	VOLUME={177},
	NUMBER={2},
	PAGES={153--172},
	YEAR={2006},
}

@book {Kigami,
	AUTHOR = {Kigami, J.},
	TITLE = {Analysis on {F}ractals},
	SERIES = {Cambridge Tracts in Mathematics},
	VOLUME = {143},
	PUBLISHER = {Cambridge University Press, Cambridge},
	YEAR = {2001},
	PAGES = {viii+226},
}

@Article{Massopust2019,
	AUTHOR = {Massopust, P.},
	TITLE = {Non-Stationary Fractal Interpolation},
	JOURNAL = {Mathematics},
	VOLUME = {7},
	YEAR = {2019},
	NUMBER = {8},
}

@article {MeirKeeler,
	AUTHOR = {Meir, A. and Keeler, E.},
	TITLE = {A theorem on contraction mappings},
	JOURNAL = {J. Math. Anal. Appl.},
	FJOURNAL = {Journal of Mathematical Analysis and Applications},
	VOLUME = {28},
	YEAR = {1969},
	PAGES = {326--329},
}

@article {Mihail,
	AUTHOR = {Mihail, A. and Miculescu, R.},
	TITLE = {Applications of fixed point theorems in the theory of
	generalized {IFS}},
	JOURNAL = {Fixed Point Theory Appl.},
	FJOURNAL = {Fixed Point Theory and Applications},
	YEAR = {2008},
	PAGES = {Art. ID 312876, 11},}

@article {MondalJha,
	AUTHOR = {Mondal, A. I. and Jha, S.},
	TITLE = {Non-stationary {$\alpha$}-fractal functions and their
	dimensions in various function spaces},
	JOURNAL = {Indag. Math. (N.S.)},
	FJOURNAL = {Koninklijke Nederlandse Akademie van Wetenschappen.
	Indagationes Mathematicae. New Series},
	VOLUME = {35},
	YEAR = {2024},
	NUMBER = {1},
	PAGES = {159--180},
}

@article {Pasupathi,
	AUTHOR = {Pasupathi, R. and Miculescu, R.},
	TITLE = {A very general framework for fractal interpolation functions},
	JOURNAL = {J. Math. Anal. Appl.},
	FJOURNAL = {Journal of Mathematical Analysis and Applications},
	VOLUME = {534},
	YEAR = {2024},
	NUMBER = {2},
	PAGES = {Paper No. 128093, 17},
}

@article {Ri2017,
	AUTHOR = {Ri, S.},
	TITLE = {A new nonlinear fractal interpolation function},
	JOURNAL = {Fractals},
	FJOURNAL = {Fractals. Complex Geometry, Patterns, and Scaling in Nature
	and Society},
	VOLUME = {25},
	YEAR = {2017},
	NUMBER = {6},
	PAGES = {1750063, 12},
}

@article {Ruan2010,
	AUTHOR = {Ruan, H. J.},
	TITLE = {Fractal interpolation functions on post critically finite
	self-similar sets},
	JOURNAL = {Fractals},
	FJOURNAL = {Fractals. Complex Geometry, Patterns, and Scaling in Nature
	and Society},
	VOLUME = {18},
	YEAR = {2010},
	NUMBER = {1},
	PAGES = {119--125},
}

@article {Ri2011,
	AUTHOR = {Ri, S. and Ruan, H. J.},
	TITLE = {Some properties of fractal interpolation functions on
	{S}ierpinski gasket},
	JOURNAL = {J. Math. Anal. Appl.},
	FJOURNAL = {Journal of Mathematical Analysis and Applications},
	VOLUME = {380},
	YEAR = {2011},
	NUMBER = {1},
	PAGES = {313--322},
}

@article {Ri2020,
	AUTHOR = {Ri, S.},
	TITLE = {Fractal functions on the {S}ierpinski gasket},
	JOURNAL = {Chaos Solitons Fractals},
	FJOURNAL = {Chaos, Solitons \& Fractals},
	VOLUME = {138},
	YEAR = {2020},
	PAGES = {110142, 10},
}

@article{Sahu,
	TITLE={On the box-counting dimension of graphs of harmonic functions on the {S}ierpiński gasket},
	AUTHOR={Sahu, A. and Priyadarshi, A.},
         JOURNAL = {J. Math. Anal. Appl.},
	FJOURNAL={Journal of Mathematical Analysis and Applications},
	VOLUME={487},
	NUMBER={2},
	PAGES={123862},
	YEAR={2020}, }

@article {Secelean,
	AUTHOR = {Secelean, N.},
	TITLE = {The fractal interpolation for countable systems of data},
	JOURNAL = {Univ. Beograd. Publ. Elektrotehn. Fak. Ser. Mat.},
	FJOURNAL = {Univerzitet u Beogradu. Publikacije Elektrotehni\v ckog
	Fakulteta. Serija Matematika},
	VOLUME = {14},
	YEAR = {2003},
	PAGES = {11--19},
}

@article {Secelean2013,
	AUTHOR = {Secelean, N.},
	TITLE = {Iterated function systems consisting of {$F$}-contractions},
	JOURNAL = {Fixed Point Theory Appl.},
	FJOURNAL = {Fixed Point Theory and Applications},
	YEAR = {2013},
	PAGES = {2013:277, 13},
}

@book{Strichartz,
	AUTHOR    = {R. S. Strichartz},
	TITLE     = {Differential Equations on Fractals: A Tutorial},
	YEAR      = {2006},
	PUBLISHER = {Princeton University Press},
	ADDRESS   = {Princeton, NJ},
}

@article{Frieberg,
	AUTHOR = {U. R. Freiberg and M. R. Lancia},
	TITLE = {Energy Form on a Closed Fractal Curve},
	JOURNAL = {Journal for Analysis and its Applications},
	FJOURNAL = {Journal for Analysis and its Applications},
	YEAR={2004},
	PAGES = {115-137},
	VOLUME={23}, 
	NUMBER={1},
	
}

@article{Liang,
	author  = {Liang, Z. and Ruan, H. J.},
	title   = {Oscillation and dimensions of harmonic functions
	on the {Sierpi\'nski} gaskets},
	journal = {Fractals},
	volume  = {33},
	number  = {3},
	pages   = {2450134},
	year    = {2025},
	doi     = {10.1142/S0218348X24501342},
}

@article {Chand1,
    AUTHOR = {Sahu, S. K. and Chand, A. K. B. and Kumar, D.},
     TITLE = {Fractal functions based on piecewise continuous functions and
              box dimension results in {H}\"older settings},
   JOURNAL = {Commun. Nonlinear Sci. Numer. Simul.},
  FJOURNAL = {Communications in Nonlinear Science and Numerical Simulation},
    VOLUME = {163},
      YEAR = {2026},
     PAGES = {Paper No. 110509, 24},}

@article {Nasim,
    AUTHOR = {Akhtar, Md. N. and G. P. Prasad, M. and Navascu\'es, M.
              A.},
     TITLE = {Box dimension of {$\alpha$}-fractal function with variable
              scaling factors in subintervals},
   JOURNAL = {Chaos Solitons Fractals},
  FJOURNAL = {Chaos, Solitons \& Fractals},
    VOLUME = {103},
      YEAR = {2017},
     PAGES = {440--449},}
\end{document}